\documentclass[a4paper,12pt]{article}

\usepackage[ansinew]{inputenc}
\usepackage{amsfonts}
\usepackage{amssymb,amsmath}
\usepackage{mathbbol}
\usepackage{latexsym}

\newtheorem{thm}{Theorem}[section]
\newtheorem{dfn}[thm]{Definition}
\newtheorem{lem}[thm]{Lemma}
\newtheorem{prop}[thm]{Proposition}

\newtheorem{rem}[thm]{Remark}

\newcommand{\RR}{\mathbb{R}}
\newcommand{\PP}{\mathbb{P}}
\newcommand{\EE}{\mathbb{E}}

\title{Limit theorems for conditioned multitype Dawson-Watanabe processes \\
and Feller diffusions}

\author{
Nicolas \textsc{ Champagnat}\\
{\footnotesize Project-team TOSCA}\\[-0.15cm]
{\footnotesize INRIA Sophia Antipolis}\\[-0.15cm]
{\footnotesize 2004 route des lucioles - BP 93}\\[-0.15cm]
{\footnotesize 06902 Sophia Antipolis, France}\\[-0.15cm]
{\footnotesize e-mail~:~Nicolas.Champagnat@sophia.inria.fr}\\
\and 
Sylvie \textsc{R\oe lly 
\footnote{On leave of absence Centre de Math\'ematiques Appliqu\'ees, UMR C.N.R.S. 7641, \'Ecole
  Polytechnique, 91128 Palaiseau C\'edex, France.}
}\\
{\footnotesize Institut f\"ur Mathematik}\\[-0.15cm]
{\footnotesize Universit\"at Potsdam}\\[-0.15cm]
{\footnotesize Am Neuen Palais 10}\\[-0.15cm]
{\footnotesize 14469 Potsdam, Germany}\\[-0.15cm]
{\footnotesize e-mail~:~roelly@math.uni-potsdam.de}\\
}

\date{\today}

\begin{document}

\maketitle

\enlargethispage*{1000pt}

\begin{abstract} 
  A multitype Dawson-Watanabe process is conditioned, in subcritical
  and critical cases, on non-extinction in the remote future. On every
  finite time interval, its distribution is absolutely continuous 
  with respect to the law of the unconditioned
  process. A martingale problem characterization is also given. 
  Several results on the long time behavior of the conditioned mass
  process---the conditioned multitype Feller branching diffusion---are
  then proved.  The general case is considered first, where the
  mutation matrix which models the interaction between the types, is
  irreducible.  Several two-type models with decomposable mutation
  matrices are also analyzed.
\end{abstract}
\bigskip

\noindent
AMS 2000 Subject Classifications~: 60J80, 60G57.\\
KEY-WORDS: multitype measure-valued branching processes; conditioned
Dawson-Watanabe process; critical and subcritical Dawson-Watanabe
process; conditioned Feller diffusion; remote survival; long time behavior.\\

\newpage

\section*{Introduction}
The paper focuses on some conditioning of the measure-valued process
called multitype Dawson-Watanabe (MDW) process, and on its mass
process, the well-known multitype Feller (MF) diffusion.  We consider
the critical and subcritical cases, in which, for any finite initial
condition, the MF diffusion vanishes in finite time, that is the MDW
process dies out a.s.  In these cases, it is interesting to condition
the processes to stay alive forever
- an event which we call {\bf remote survival}, see the exact definition in (\ref{eq:defP*}). \\
Such a study was initiated for the monotype Dawson-Watanabe process by
A. Rouault and the second author in \cite{RR89CRAS} (see also
\cite{EvPe}, \cite{Ev} and \cite{EthWi} for the study of various
aspects of conditioned monotype superprocesses). Their results were a
generalization at the level of measure-valued processes 
of the pioneer work of Lamperti and Ney (\cite{LN68TPA}, Section 2), 
who studied the same questions applied to Galton-Watson processes.\\
We are interested here in the {\bf multitype} setting which is much different from 
the monotype one.  The mutation matrix $D$ introduced in
(\ref{eq:PDE-U}),
which measures the quantitative interaction between types, will play a crucial role.\\
We now briefly describe the contents of the paper.  The model is
precisely defined in the first section.  In the second section we
define the conditioned MDW process, express its law as a locally absolutely
continuous measure with respect to the law of the unconditioned
process, write explicitly the martingale problem it satisfies and give
the form of its Laplace functional; all this in the case of an
irreducible mutation matrix.  Since $D$ is irreducible, all the types
communicate and conditioning by remote survival is equivalent to
conditioning by the non-extinction of only one type (see Remark
\ref{rem:other-cond}).  The third section is devoted to the long time
behavior of the mass of the conditioned  MDW process, which is then a conditioned MF diffusion.
  First the monotype case is analyzed (it
was not considered in \cite{RR89CRAS}), and then the irreducible
multitype case.  We also prove that both limits interchange: the long
time limit and the conditioning by long time survival (see
Theorem~\ref{thm:interv-multi}).  In the last section we treat the
same questions as in Section~\ref{sec:long-time}
for various reducible 2-types models. Since $D$ is decomposable, the
two types can have very different behaviors, that also depend on the
precise conditioning that is considered (see Section \ref{subsec:first model}).

\section{The model}
\label{sec:intro}

In this paper, we will assume for simplicity that the (physical) space
is $\RR$.  $k$ is the number of types.  Any $k$-dimensional vector $u
\in \RR^k$ is denoted by $(u_1;\cdots;u_k)$.  $\bf 1$ will denote the
vector $(1;\ldots;1)\in\mathbb{R}^k$.  $\|u\|$ is the euclidean norm
of $u \in \RR^k$ and $(u, v)$ the scalar product between
$u$ and $v$ in $\RR^k$.
If $u \in \RR^k$, $|u|$ is the vector in $\RR^k$ with coordinates
$|u_i|, 1\leq i\leq k$.\\
We will use the notations $u>v$ (resp.\ $u\geq v$) when $u$ and $v$
are vectors or matrices such that $u-v$ has positive (resp.\
non-negative) entries.\\
Let $C_b(\RR,\RR^k)$ denote the space of $\RR^k$-valued continuous
bounded functions on $\RR$.  By $C_b(\RR,\RR^k)_+$ we denote the set
of non-negative elements of
$C_b(\RR,\RR^k)$.\\
$M(\RR)$ is the set of finite positive measures on $\RR$,
and $M(\RR)^k$ the set of $k$-dimensional vectors of finite positive measures.\\
The duality between measures and functions will be denoted by
$\langle\cdot,\cdot \rangle$~: $\langle \nu, f \rangle := \int fd\nu$
if $\nu \in M(\RR)$ and $f$ is defined on $\RR$, and in the vectorial
case
$$
\langle (\nu_1;\ldots;\nu_k),(f_1;\ldots;f_k)\rangle := \sum_{i=1}^k\int f_id\nu_i 
 = \big( (\langle \nu_1,f_1 \rangle;\ldots;\langle \nu_k,f_k \rangle),{\bf 1}\big) 
$$ 
for $\nu=(\nu_1;\ldots;\nu_k) \in M(\RR)^k$ and $f=(f_1;\ldots;f_k)
\in C_b(\RR,\RR^k)$. For any $\lambda\in\RR^k$, the constant function
of $C_b(\RR,\RR^k)$ equal to $\lambda$ will be also denoted by
$\lambda$.
\bigskip

A multitype Dawson-Watanabe process with \emph{mutation matrix}
$D=(d_{ij})_{1\leq i,j\leq k}$ is a continuous $M(\RR)^k$-valued
Markov process whose law $\PP$ on the canonical space
$(\Omega:=C(\RR_+, M(\RR)^k),(X_t)_{t \geq 0}, ({\cal F}_t)_{t \geq 0})$ has as transition
Laplace functional
\begin{equation}
  \label{eq:LaplaceMDW}
  \forall f\in C_b(\RR,\RR^k)_+, \quad \EE(\exp -\langle
  X_t,f\rangle\mid X_0=m)=\exp-\langle m, U_tf\rangle
\end{equation}
where $U_t f \in C_b(\RR,\RR^k)_+$, the so-called cumulant semigroup,
is the unique solution of the non-linear PDE
\begin{equation}
  \label{eq:PDE-U}
  \left\{
    \begin{array}{rl}
      \displaystyle   \frac{\partial (U_tf)}{\partial t} & = 
      \displaystyle \Delta  U_tf+DU_tf- \frac{c}{2}(U_tf)^{\odot 2}\\
      U_0f & =  f .
    \end{array}
  \right.
\end{equation}
Here, $u\odot v$ denotes the componentwise product $(u_iv_i)_{1\leq
  i\leq k}$ of two $k$-dimensional vectors $u$ and $v$ and $u^{\odot
  2}=u\odot u$.  To avoid heavy notation, when no confusion is
possible, we do not write differently column and row vectors when
multiplied by a matrix. In particular, in the previous
equation, $Du$ actually stands for $Du'$.\\
The MDW process arises as the diffusion limit of a sequence of
particle systems $(\frac{1}{K}N^K)_K$, where $N^K$ is an appropriate
rescaled multitype branching Brownian particle system (see e.g.\
\cite{GL-M90AAP} and \cite{GR91SPA}, or \cite{Wa68JMKU} for the
monotype model): after an exponential lifetime with parameter $K$,
each Brownian particle splits or dies, in such a way that the number
of offsprings of type $j$ produced by a particle of type $i$ has as
(nonnegative) mean \mbox{$\delta_{ij}+ \frac{1}{K} d_{ij}$} and as second factorial
moment $c$ ($\delta_{ij}$ denotes the Kronecker function, equal to 1
if $i=j$ and to 0 otherwise).  Therefore, the average number of
offsprings of each particle is asymptotically one and the matrix $D$
measures the (rescaled) discrepancy between the mean matrix and the
identity matrix $I$,
which corresponds to the pure critical case of independent types. \\
For general literature on DW processes we refer the reader e.g.\ to
the lectures of D. Dawson \cite{Da93StFlour} and E. Perkins
\cite{Pe02StFlour}
and the monographs \cite{Dy94AMS} and~\cite{Eth00AMS}.\\
Let us remark that we introduced a variance parameter $c$ which is
type-independent.  In fact we could replace it by a vector
$c=(c_1;\cdots ; c_k)$, where $c_i$ corresponds to type $i$.  If
$\inf_{1\leq i\leq k} c_i>0$, then all the results of this paper are
still true.  We decided to take $c$ independent of the type to
simplify the notation.

When the mutation matrix $D=(d_{ij})_{1\leq i,j\leq k}$ is not diagonal, it represents the
interaction between the types, which justifies its name. 
Its non diagonal elements are non-negative. These
matrices are sometimes called Metzler-Leontief matrices in
financial mathematics (see \cite{Sen} § 2.3 and the bibliography
therein).  Since there exists a positive constant
$\alpha$ such that $D+\alpha I\geq 0$, it follows from
Perron-Frobenius theory that $D$ has a real eigenvalue $\mu$ such that no
other eigenvalue of $D$ has its real part exceeding $\mu$.  Moreover,
the matrix $D$ has a non negative right eigenvector associated to the
eigenvalue $\mu$ (see e.g.~\cite{Ga86}, Satz 3 § 13.3 or \cite{Sen}
Exercise 2.11).  The cases $\mu<0$, $\mu=0$ and $\mu>0$ correspond
respectively to a \emph{subcritical}, \emph{critical} and
\emph{supercritical} processes. \\
In the present paper, we only consider the case $\mu \leq 0$, in which
the MDW dies out a.s.\ (see Jirina~\cite{Ji64Prague}).

\section{(Sub)critical irreducible MDW process conditioned by remote survival.}
\label{sec:irred}

Let us recall the definition of irreducibility of a matrix.
\begin{dfn}
  \label{def:irredetprimit}     
  A square matrix $D$ is called \emph{irreducible} if there is no
  permutation matrix $Q$ such that $Q^{-1}DQ$ is block triangular.
\end{dfn}
In all this section and in the next one, the mutation matrix $D$ is
assumed to be irreducible.  By Perron-Frobenius' theorem (see
e.g.~\cite{Sen} Theorem 1.5 or~\cite{Ga86}, Satz 2 §13.2, based on
\cite{Pe1907} and \cite{Fro1912}), the eigenspace associated to the
maximal real eigenvalue $\mu$ of $D$ is one-dimensional.  We will
always denote its generating right (resp.\ left) eigenvector by $\xi$
(resp.\ by $\eta$) with the normalization conventions $(\xi,{\bf
  1})=1$ and $(\xi,\eta)=1$.
All the coordinates of both vectors $\xi$ and $\eta$ are positive.

\subsection{The conditioned process as a $h$-process}
\label{sec:cond}

The natural way to define the law $\PP^*$ of the MDW process
conditioned to never die out is by
\begin{equation}
  \label{eq:defP*}
   \forall B\in{\cal F}_t, \quad
  \PP^*(B) := \lim_{\theta\rightarrow\infty} \PP(B\mid\langle
  X_{t+\theta},{\bf 1} \rangle>0)
\end{equation}
if this limit exists.

The following Theorem~\ref{thm:cond} proves that $\PP^*$ is
well-defined by~(\ref{eq:defP*}) and is a probability measure on
${\cal F}_t$ absolutely continuous with respect to $\PP$. Furthermore,
the density is a martingale, so that $\PP^*$ can be extended to
$\vee_{t\geq 0} {\cal F}_t$, defining a Doob $h$-transform of $\PP$
(see the seminal work~\cite{Mey} on $h$-transforms and~\cite{Ov} for
applications to monotype DW processes).

\begin{thm}
  \label{thm:cond}
  Let $\PP$ be the distribution of a critical or subcritical MDW
  process characterized by \textup{(\ref{eq:LaplaceMDW})}, with an
  irreducible mutation matrix $D$ and initial measure $m \in M(\RR)^k
  \setminus \{0\}$. Then, the limit in~\textup{(\ref{eq:defP*})}
  exists and defines a probability measure $\PP^*$ on $\vee_{t\geq 0}
  {\cal F}_t$ such that, for any $t>0$,
  \begin{equation}
    \label{eq:cond}
    \PP^*\bigr|_{{\cal F}_t} = \frac{\langle X_t,\xi\rangle}{\langle m,\xi\rangle} \, 
    e^{-\mu t} \quad  \PP\bigr|_{{\cal F}_t} 
\end{equation}
where $\xi \in \RR^k$ is the unitary right eigenvector associated to
the maximal real eigenvalue $\mu$ of $D$.
\end{thm}

\paragraph{Proof of Theorem~\ref{thm:cond}}
By definition, for $B\in{\cal F}_t$,
\begin{equation*}
  \EE(\mathbb{1}_B\mid\langle X_{t+\theta},{\bf 1} \rangle>0) 
  = \frac{\EE\big( \mathbb{1}_B(1-\PP(\langle X_{t+\theta},{\bf 1}
    \rangle=0\mid{\cal F}_t))\big)}
  {1-\PP(\langle X_{t+\theta},{\bf 1} \rangle=0)} .
\end{equation*}
For any time $s>0$, $x_s:=(\langle X_{s,1}, 1\rangle; \dots; \langle
X_{s,k}, 1\rangle) $, the total mass at time $s$ of the MDW
process---a multitype Feller diffusion---is a continuous
$\RR_+^k$-valued process with initial value $x=(\langle m_1, 1\rangle;
\dots ; \langle m_k, 1\rangle)$ characterized by its transition
Laplace transform
\begin{equation}
  \label{eq:TL-x}
  \displaystyle 
  \forall\lambda\in\RR_+^k,\quad 
  \EE( e^{ - (x_t, \lambda)} \mid x_0=x)=e^{- (x,u_t^\lambda)} .
\end{equation}
Here, $u_t^\lambda=(u_{t,1}^\lambda;\ldots
;u_{t,k}^\lambda):=U_t\lambda$ satisfies the non-linear differential
system
\begin{equation}
  \label{eq:ODE-u}
  \left\{
    \begin{array}{rl}
      \displaystyle  \frac{du_t^\lambda}{dt} &
      =Du_t^\lambda-\displaystyle \frac{c}{2}(u_t^\lambda)^{\odot 2}\\
      u_0^\lambda & =\lambda ,
    \end{array}
  \right.
\end{equation}
or componentwise
\begin{equation*}
  \forall i\in\{1,\ldots,k\},\quad
  \frac{du_{t,i}^\lambda}{dt}=\sum_{j=1}^kd_{ij}u_{t,j}^\lambda
  -\frac{c}{2}(u_{t,i}^\lambda)^2,\quad u_{0,i}^\lambda=\lambda_i .
\end{equation*}
Then,
$$
\displaystyle 
\PP(\langle X_{s},{\bf 1} \rangle=0) = 
\lim_{\lambda\hookrightarrow\infty} \EE(e^{-\langle X_s,\lambda \rangle} \mid X_0=m)
= e^{-(x,\lim_{\lambda\hookrightarrow\infty}u_s^\lambda)}
$$
where $\lambda \hookrightarrow \infty$ means that all coordinates of
$\lambda$ go to $+\infty$.  Using the Markov property of the MDW
process, one obtains
\begin{equation} 
  \label{eq:RR-1}
  \EE(\mathbb{1}_B\mid\langle X_{t+\theta},{\bf 1} \rangle>0) 
  =\frac{\EE\Big( \mathbb{1}_B \big(
    1-e^{-(x_t,\lim_{\lambda\hookrightarrow\infty}u_\theta^\lambda)}
    \big) \Big)}
  {1-e^{-(x,\lim_{\lambda\hookrightarrow\infty}u_{t+\theta}^\lambda)}} .
\end{equation}
\medskip

In the monotype case ($k=1$), $u_t^\lambda$ can be computed
explicitly (see Section \ref{sec:long-time-monotype}), but this is
not possible in the multitype case. Nevertheless, one can obtain upper
and lower bounds for $u_t^\lambda$.  This is the goal of the
following two lemmas, the proofs of which are postponed after the end of the
proof of Theorem~\ref{thm:cond}.
\begin{lem}
  \label{lem:pte-u}
  Let $u_t^\lambda=(u_{t,1}^\lambda;\ldots ;u_{t,k}^\lambda)$ be the
  solution of \textup{(\ref{eq:ODE-u})}.
  \begin{description}
  \item[\textmd{(i)}] For any $\lambda\in\RR^k_+\setminus\{0\}$ and
    any $t>0$, $u_t^\lambda>0$.
  \item[\textmd{(ii)}] Let $\displaystyle
    C_t^\lambda:=\sup_{1\leq i\leq k}\frac{u_{t,i}^\lambda}{\xi_i}$
    and   $\underline{\xi}:=\inf_i\xi_i$. For $t> 0$ and $ \lambda\in \RR^k_+$,\\
    - in the critical case ($\mu=0$)
    \begin{equation}
      \label{eq:C_tcritique}
      C^\lambda_t\leq  \frac{C_0^\lambda}{1+\frac{c
          \underline\xi}{2}C_0^\lambda\:t} \quad
      \textrm{and therefore}   
      \quad  \sup_{\lambda\in \RR^k_+}C^\lambda_t\leq
      \frac{2}{c\underline\xi\:t}  
    \end{equation}
    - in the subcritical case ($\mu<0$)
    \begin{equation}
      \label{eq:C_tsubcritique}
      C^\lambda_t\leq    \frac{C_0^\lambda e^{\mu
          t}}{1+\frac{c\underline\xi}{2|\mu|}C_0^\lambda(1-e^{\mu t})} \quad
      \textrm{and therefore} \quad
      \sup_{\lambda\in \RR^k_+}C^\lambda_t\leq 
      \frac{2|\mu|e^{\mu t}}{c\underline\xi(1-e^{\mu t})} 
    \end{equation}   
  \item[\textmd{(iii)}] Let $\displaystyle
    B_t^\lambda:=\inf_{1\leq i\leq k}\frac{u_{t,i}^\lambda}{\xi_i}$
    and $\bar\xi:=\sup_i\xi_i$.  Then
    \begin{equation}
      \label{eq:B_t}
      \forall t\geq 0,\ \lambda\in \RR^k_+,\quad B^\lambda_t\geq\left\{
        \begin{array}{ll}
          \displaystyle \frac{B_0^\lambda}{1+\frac{c\bar\xi}{2}
            B_0^\lambda\:t} & \mbox{if\ 
          }\mu=0 \\
          \displaystyle \frac{B_0^\lambda e^{\mu
              t}}{1+\frac{c\bar\xi}{2|\mu|} B_0^\lambda (1-e^{\mu t})} & \mbox{if\
          }\mu<0.
        \end{array}
      \right.
    \end{equation}
  \item[\textmd{(iv)}] For any $\lambda\in \RR^k_+$ and $t\geq 0$,
    \begin{equation}
      \label{eq:LB-u}
      u_t^\lambda\geq\left\{
        \begin{array}{ll}
          \left(1+\frac{c \underline\xi}{2}C_0^\lambda
            \:t\right)^{-\bar\xi/\underline\xi}e^{Dt}\lambda
          & \mbox{if\ }\mu=0 \\
          \left(1+\frac{c\underline\xi}{2|\mu|}C_0^\lambda(1-e^{\mu
              t})\right)^{-\bar\xi/\underline\xi}e^{Dt}\lambda & 
          \mbox{if\ }\mu<0.   
        \end{array}
      \right.
    \end{equation}
  \end{description}
\end{lem}

The main difficulty in the multitype setting comes from the non-commutativity of
matrices. For example~(\ref{eq:ODE-u}) can be expressed as
$\frac{du^\lambda_t}{dt}=(D+A_t)u^\lambda_t$ where the matrix $A_t$ is
diagonal with $i$-th diagonal element $cu^\lambda_{t,i}/2$. However,
since $D$ and $A_t$ do not commute, it is not possible to express
$u^\lambda_t$ in terms of the exponential of $\int_0^t (D+A_s)\, ds$.
 The following lemma gives the main tool we use to solve
this difficulty.
\begin{lem}
  \label{lem:diff-ineq}
  Assume that $t \mapsto f(t)\in\mathbb{R}$ is a continuous function
  on $\RR_+$ and $t \mapsto u_t\in \RR^k$ is a differentiable
  function on $\RR_+$. Then
  $$
  \frac{du_t}{dt}\geq (D + f(t) I) u_t, \quad \forall t\geq 0
  \quad \Longrightarrow \quad
  u_t\geq\exp\Big(\int_0^t (D + f(s) I)\, ds \Big) u_0, \quad \forall t\geq 0
  $$
\end{lem}
\medskip

For any $ 1\leq i\leq k $, applying~(\ref{eq:TL-x}) with $x=e^i$ where
$e_j^i=\delta_{ij},1\leq j\leq k $, one easily deduces
the existence of a limit in $[0,\infty]$ of $u_{t,i}^\lambda$ when
$\lambda \hookrightarrow \infty$.
Moreover, by Lemma~\ref{lem:pte-u}~(ii) and~(iii), for any
$t>0$,
$$
0<\frac{2 f(\theta)}{c\bar{\xi}}\leq
\lim_{\lambda\hookrightarrow\infty}u^\lambda_\theta\leq
\frac{2 f(\theta)}{c\underline{\xi}}<+\infty
$$
where $f(\theta)=1/\theta$ if $\mu=0$ or
$f(\theta)=|\mu|e^{\mu\theta}/(1-e^{\mu\theta})$ if $\mu<0$. Therefore
$\lim_{\theta\rightarrow\infty}\lim_{\lambda\hookrightarrow\infty}
u_\theta^\lambda=0$ and, for sufficiently large $\theta$,
$$
\frac{1-e^{-(x_t,\lim_{\lambda\hookrightarrow\infty}u_\theta^\lambda)}}
{1-e^{-(x,\lim_{\lambda\hookrightarrow\infty}u_{t+\theta}^\lambda)}}\leq
K\frac{(x_t,\mathbf{1})}{(x,\mathbf{1})}
$$
for some constant $K$ that may depend on $t$ but is independent of
$\theta$. Since $\EE\langle X_t,\bf 1\rangle<\infty$ for any $t\geq 0$
(see \cite{GL-M90AAP} or \cite{GR91SPA}), Lebesgue's dominated
convergence theorem can be applied to make a first-order expansion in
$\theta$ in~(\ref{eq:RR-1}). This yields that the density with respect
to $\PP$ of $\PP$ conditioned on the non-extinction at time $t+\theta$
on ${\cal F}_t$, converges in $L^1(\PP)$ when
$\theta\rightarrow\infty$ to
\begin{equation}
  \label{eq:bla}
   \big( x_t,\lim_{\theta\rightarrow\infty}
   \frac{\lim_{\lambda\hookrightarrow\infty}u_\theta^\lambda}
   { (x,\lim_{\lambda\hookrightarrow\infty}u_{t+\theta}^\lambda)} \big)
\end{equation}
if this limit exists.\\
We will actually prove that
\begin{equation}
  \label{eq:limitunifu}
  \lim_{\theta\rightarrow\infty}  \sup_{\lambda\not= 0}
  \|\frac{1}{(x,u_{t+\theta}^\lambda)} u_\theta^\lambda- 
  \frac{e^{-\mu t}}{(x,\xi)} \xi \| =0.
\end{equation}
This will imply that the limits in $\theta$ and in $\lambda$ can be
exchanged in~(\ref{eq:bla}) and thus
$$
\lim_{\theta\rightarrow\infty} \EE(\mathbb{1}_B\mid\langle X_{t+\theta},{\bf 1} \rangle>0) 
=  e^{-\mu t} \,  \EE\Big( \mathbb{1}_B \frac{(x_t ,\xi)}{(x ,\xi)}\Big) = 
e^{-\mu t} \,  \EE\Big( \mathbb{1}_B \frac{\langle X_t ,\xi\rangle}{\langle m
  ,\xi\rangle}\Big),
$$
completing the proof of Theorem~\ref{thm:cond}.

\noindent \underline{Subcritical case: $\mu<0$}\\
As a preliminary result, observe that, since $D$ has
nonnegative nondiagonal entries, there exists $\alpha>0$ such that
$D+\alpha I\geq 0$, and then $\exp(Dt)\geq 0$.

Since $\displaystyle \frac{du_t^\lambda}{dt}\leq Du_t^\lambda$, we
first remark by Lemma~\ref{lem:diff-ineq} (applied to $-u^\lambda_t$),
that
$$
\forall t\geq 0,\quad u_t^\lambda\leq e^{Dt}\lambda.
$$ 
Second, it follows from Lemma~\ref{lem:pte-u} (iv) that
\begin{align*}
  e^{Dt}\lambda - u_t^\lambda &
  \leq\bigg(1-\big(1+\frac{c\underline\xi}{2|\mu|}
  C_0^\lambda(1-e^{\mu t})\big)^{-\bar\xi/\underline\xi}
  \bigg)e^{Dt}\lambda \\
  & \leq\frac{c\bar\xi}{2|\mu|}C_0^\lambda(1-e^{\mu t}) \, e^{Dt}\lambda .
\end{align*}
Therefore, since $C_0^\lambda=\sup_i\lambda_i/\xi_i$, there exists a
constant $K$ independent of $\lambda$ such that
\begin{equation}
  \label{eq:equiv}
  \forall\lambda\geq 0,\ \forall t\geq 0,\quad
  e^{Dt}\lambda - u_t^\lambda\leq K\|\lambda\|e^{Dt}\lambda.
\end{equation}
In particular,
$$
\|\lambda\|\leq \frac{1}{2K} \quad\Rightarrow\quad
u^\lambda_t\geq\frac{1}{2}e^{Dt}\lambda\quad\forall t\geq 0.
$$
Third, it follows from Lemma~\ref{lem:pte-u} (ii) that there exists
$t_0$ such that
\begin{equation*}
  \label{eq:def-t_0}
  \forall t\geq t_0,\ \forall\lambda\geq 0,\quad \|u_t^\lambda\|\leq \frac{1}{2K}.
\end{equation*}
Fourth, as a consequence of Perron-Frobenius' theorem, the exponential
matrix $e^{Dt}$ decreases like $e^{\mu t}$ for $t$ large in the
following sense: as $t\rightarrow \infty$,
\begin{equation} 
  \label{equivalentexpDt}
  \exists \gamma>0, \quad  e^{Dt}=
  e^{\mu t}P+O(e^{(\mu-\gamma)t})
\end{equation}
where $P:=(\xi_i\eta_j)_{1\leq i,j\leq k}$
(see \cite{Sen} Theorem 2.7). Therefore, there exists $\theta_0$ such that
\begin{equation*}
  \forall t \geq\theta_0,
  \quad \frac{1}{2} e^{\mu t} P \leq e^{Dt}\leq  2e^{\mu t}P.
\end{equation*}
Last, there exists a positive constant $K'$ such that
\begin{equation*}
  \forall u,v\in\mathbb{R}^k_+,\quad  (v,Pu) = (u,\eta)(v,\xi)\geq
  \underline\xi \, \underline\eta \, (u,{\bf 1})
  \,   (v ,{\bf 1})\, \geq K'\|u\|\|v\| .
\end{equation*}
Combining all the above inequalities, we get for any $a\in\RR^k_+$,
$b, \lambda \in\RR^k_+\setminus\{0\}$
and for any $\theta\geq\theta_0$,
\begin{align}
  & \left|\frac{( a,u_{t_0+\theta}^\lambda)}
    {( b,u_{t_0+\theta+t}^\lambda)}
    -\frac{( a,e^{D\theta}u_{t_0}^\lambda)}
    {( b,e^{D(\theta+t)}u_{t_0}^\lambda)}\right| \notag \\
  & \qquad\qquad \leq
  \frac{( a,|u_{t_0+\theta}^\lambda-e^{D\theta}u_{t_0}^\lambda|)}
  {( b,u_{t_0+\theta+t}^\lambda)}
  +\frac{( a,e^{D\theta}u_{t_0}^\lambda)
    ( b,|u_{t_0+\theta+t}^\lambda-e^{D(\theta+t)}u_{t_0}^\lambda|)}
  {( b,u_{t_0+\theta+t}^\lambda)
    ( b,e^{D(\theta+t)}u_{t_0}^\lambda)} \notag \\
  & \qquad\qquad \leq
  \frac{2K\|a\|\|u_{t_0}^\lambda\|\|e^{D\theta}u_{t_0}^\lambda\|}
  {( b,e^{D(\theta+t)}u_{t_0}^\lambda)}
  +\frac{2K\|a\|\|e^{D\theta}u_{t_0}^\lambda\|
    \|b\|\|u_{t_0}^\lambda\|\|e^{D(\theta+t)}u_{t_0}^\lambda\|}
  {( b,e^{D(\theta+t)}u_{t_0}^\lambda)^2} \notag \\
  & \qquad\qquad \leq \bar K \|a\|\|u_{t_0}^\lambda\|e^{-\mu t}\left(
    \frac{\|Pu_{t_0}^\lambda\|}{( b,Pu_{t_0}^\lambda)}
    +\frac{\|b\|\|Pu_{t_0}^\lambda\|^2}{( b,Pu_{t_0}^\lambda)^2
    }\right) \notag \\
  & \qquad\qquad \leq \bar K e^{-\mu
    t}\frac{\|a\|}{\|b\|}\|u_{t_0}^\lambda\| \label{eq:calc}
\end{align}
where the constants $\bar K$ may vary from line to line, but are
independent of $\lambda$ and $t_0$.

Now, let $t_0(\theta)$ be an increasing function of $\theta$ larger
than $t_0$ such that $t_0(\theta)\rightarrow \infty$ when
$\theta\rightarrow\infty$. By Lemma~\ref{lem:pte-u} (ii),
$\|u_{t_0(\theta)}^\lambda\|\rightarrow 0$ when
$\theta\rightarrow\infty$, uniformly in $\lambda\geq 0$. Then,
by~(\ref{eq:calc}), uniformly in $\lambda\geq 0$,
\begin{align*}
  \lim_{\theta\rightarrow\infty}
  \frac{( a,u_{t_0(\theta)+\theta}^\lambda)}
  {( b,u_{t_0(\theta)+\theta+t}^\lambda)}
  & =\lim_{\theta\rightarrow\infty}
  \frac{( a,e^{D\theta}u_{t_0(\theta)}^\lambda)}
  {( b,e^{D(\theta+t)}u_{t_0(\theta)}^\lambda)} \\
  & =\lim_{\theta\rightarrow\infty}
  \frac{( a,e^{\mu\theta}Pu_{t_0(\theta)}^\lambda)}
  {( b,e^{\mu(\theta+t)}Pu_{t_0(\theta)}^\lambda)} \\
  & =\lim_{\theta\rightarrow\infty}e^{-\mu t}
  \frac{(\eta,u_{t_0(\theta)}^\lambda)( a,\xi)}
  {(\eta,u_{t_0(\theta)}^\lambda)( b,\xi)} \\
  & =e^{-\mu t}\frac{( a,\xi)}{( b,\xi)}
\end{align*}
which completes the proof of Theorem~\ref{thm:cond} in the case
$\mu<0$.

\noindent \underline{Critical case: $\mu=0$}\\
The above computation has to be slightly modified.
Inequality~(\ref{eq:equiv}) becomes
\begin{align}
  |u_t^\lambda-e^{Dt}\lambda|
  \leq & \bigg(1-\big( 1+\frac{c\underline\xi}{2}
     C_0^\lambda \: t \big)^{-\bar\xi/\underline\xi}
  \bigg)e^{Dt}\lambda \notag \\
  \leq & K\|\lambda\| te^{Dt}\lambda. \label{ineq:approxu-t}
\end{align}
Therefore, the right-hand side of~(\ref{eq:calc}) has to be replaced
by
\begin{equation}
  \label{eq:rhs-calc}
  K \frac{\|a\|}{\|b\|}\|u_{t_0}^\lambda\|(\theta+t).
\end{equation}
Now, using Lemma~\ref{lem:pte-u} (iii) again, it suffices to choose a
function $t_0(\theta)$ in such a way that
$\lim_{\theta\rightarrow\infty}\theta\sup_{\lambda\geq
  0}\|u_{t_0(\theta)}^\lambda\|= 0$.
One can now complete the proof of Theorem~\ref{thm:cond} as above.\hfill$\Box$\bigskip

\paragraph{Proof of Lemma~\ref{lem:pte-u}} \hspace{10cm}

\textbf{(i)} First, observe that, by~(\ref{eq:TL-x}), $u^\lambda_t\geq
0$ for any $t\geq 0$. Next, since $D$ is nonnegative outside the
diagonal,
\begin{equation}
  \label{eq:pf-L-2.3}
  \frac{du^\lambda_{t,i}}{dt}=\sum_{j=1}^{k}d_{ij}u^\lambda_{t,j}
  -\frac{c}{2}(u^\lambda_{t,i})^2\geq (d_{ii}-\frac{c}{2}u^\lambda_{t,i})u^\lambda_{t,i}.
\end{equation}
Therefore, for any $i$ such that $\lambda_i>0$, $u^\lambda_{t,i}>0$
for any $t\geq 0$. \\
Let $I:=\{i:\lambda_i>0\}$ and
$J:=\{j:\lambda_j=0\}$. By the irreducibility of the matrix $D$, there
exist $i\in I$ and $j\in J$ such that $d_{ji}>0$. Therefore, for
sufficiently small $t>0$,
$$
\frac{du^\lambda_{t,j}}{dt}
= \sum_{l=1}^{k}d_{jl}u^\lambda_{t,l}
  -\frac{c}{2}(u^\lambda_{t,j})^2 > \frac{d_{ji}}{2}u^\lambda_{t,i}
$$
and thus $u^\lambda_{t,j}>0$ for $t>0$ in a neighborhood of
0. Moreover, as long as $u^\lambda_{t,i}>0$, for the same reason,
$u^\lambda_{t,j}$ cannot reach 0.\\
Defining $I'=I\cup\{j\}$ and $J'=J\setminus\{j\}$, there exists $i'\in
I'$ and $j'\in J'$ such that $d_{j'i'}>0$. For sufficiently small
$\varepsilon>0$, $u^\lambda_{\varepsilon,i'}>0$ and the previous
argument shows that $u^\lambda_{\varepsilon+t,j'}>0$ for $t>0$ as long
as $u^\lambda_{\varepsilon+t,i'}>0$. Letting $\varepsilon$ go to 0
yields that $u^\lambda_{t,j'}>0$ for sufficiently small $t>0$.\\
Applying the same argument inductively shows that $u^\lambda_t>0$ for
$t>0$ in a neighborhood of 0.  Using~(\ref{eq:pf-L-2.3}) again, this
property can be extended to all $t>0$.

\textbf{(ii) and (iii)} As the supremum of finitely many continuously
differentiable functions, $t\mapsto
C_t^\lambda$ is differentiable except at at most countably many
points. Indeed, it is not differentiable at time $t$ if and only if
there exist two types $i$ and $j$ such that
$u^\lambda_{t,i}/\xi_i=u^\lambda_{t,j}/\xi_j$ \textbf{and}
$d(u^\lambda_{t,i}/\xi_i)/dt\not=d(u^\lambda_{t,j}/\xi_j)/dt$. For
fixed $i$ and $j$, such points are necessarily isolated, and hence are
at most denumerable.\\
Fix a time $t$ at which $C_t^\lambda$ is differentiable and fix $i$
such that $u_{t,i}^\lambda=C_t^\lambda\xi_i$. Then
\begin{align*}
  \frac{dC_t^\lambda}{dt}\xi_i=\frac{du_{t,i}^\lambda}{dt} &
  =\sum_{j=1}^k d_{ij}u^\lambda_{t,j}-\frac{c}{2}(u_{t,i}^\lambda)^2 \\
  & \leq C^\lambda_t\sum_{j\not=i}d_{ij}\xi_j+d_{ii}u^\lambda_{t,i}
  -\frac{c}{2}(u_{t,i}^\lambda)^2 \\ &
  =C^\lambda_t(D\xi)_i-\frac{c}{2}(u_{t,i}^\lambda)^2=\mu
  C^\lambda_t\xi_i-\frac{c}{2}\xi_i^2(C_t^\lambda)^2
\end{align*}
where the inequality comes from the fact that $D$ is nonnegative
outside of the diagonal and where the third line comes from the
specific choice of the subscript $i$. Therefore,
\begin{equation} \label{eq:Ct}
  \frac{dC^\lambda_t}{dt}\leq\mu
  C^\lambda_t-\frac{c}{2}\underline\xi\:(C_t^\lambda)^2.
\end{equation}\\
Assume $\mu=0$.\\
By Point~(i), if $\lambda\not=0$, $C_t^\lambda>0$ for any $t\geq 0$ (the
case $\lambda=0$ is trivial). Then, for any $t\geq 0$, except at at
most countably many points,
$$
- \frac{dC^\lambda_t/dt}{{(C_t^\lambda})^2} \geq  \frac{c}{2} \, \underline\xi.
$$
Integrating this inequality between 0 and $t$, we get
$$
\frac{1}{C_t^\lambda} \geq \frac{1}{C_0^\lambda} + \frac{c}{2}\underline\xi t
\quad \Rightarrow \quad
C_t^\lambda \leq \frac{C_0^\lambda}{1 + \frac{c C_0^\lambda}{2}\underline\xi t}.
$$
The proof of the case $\mu<0$ can be done by the same argument
applied to $t \mapsto e^{-\mu t}C_t^\lambda$.  Inequalities (iii) are
obtained in a similar way too.

\textbf{(iv)} By definition of $C_t^\lambda$,~(\ref{eq:ODE-u}) implies
that
\begin{equation*}
  \frac{d u_t^\lambda}{dt}\geq\left(D-\frac{c}{2}\bar\xi \:C_t^\lambda
    I\right)u_t^\lambda.
\end{equation*}
Then, (iv) follows from (ii) and Lemma~\ref{lem:diff-ineq}.
\hfill$\Box$\bigskip

\paragraph{Proof of Lemma~\ref{lem:diff-ineq}}
Fix $\varepsilon>0$ and let
\begin{equation*}
  u_t^{(\varepsilon)}:=\exp\big(\int_0^t (D+f(s)I)ds\big)(u_0-\varepsilon).
\end{equation*}
Then
$$
\frac{du_t}{dt}-\frac{du_t^{(\varepsilon)}}{dt} \geq (D+f(t)I)(u_t-u_t^{(\varepsilon)}).
$$
Let $t_0:=\inf\{t\geq 0:\exists i \in\{1,\ldots,k \},
u_{t,i}<u_{t,i}^{(\varepsilon)}\}$. For any $t\leq t_0$, since $D$ is
nonnegative outside of the diagonal,
$$
\forall i\in\{1,\ldots,k \},\quad
\frac{d}{dt}(u_{t,i}-u_{t,i}^{(\varepsilon)})
\geq(d_{ii}+f(t))(u_{t,i}-u_{t,i}^{(\varepsilon)}).
$$
Since $u_0>u_0^{(\varepsilon)}$, this implies that
$u_{t}-u^{(\varepsilon)}_{t}>0$ for any $t\leq t_0$ and thus
$t_0=+\infty$.  Letting $\varepsilon$ go to 0 completes the proof of
Lemma~\ref{lem:diff-ineq}.\hfill$\Box$\bigskip

\begin{rem} \label{rem:other-cond} \textup{Since the limit in
    (\ref{eq:limitunifu}) is uniform in $\lambda$, one can choose in
    particular $\lambda=\lambda^i, 1\leq i\leq k $, where
    $\lambda^i_j=0$ for $j\not=i$.  Thus, for each type $i$,
    \begin{equation*}
      \lim_{\theta\rightarrow\infty}\frac{
        \lim_{\lambda^i_i \rightarrow\infty}u_\theta^{\lambda^i}}
      {(x, \lim_{\lambda^i_i \rightarrow\infty}u_{t+\theta}^{\lambda^i})}
      =\frac{e^{-\mu t}}{(x,\xi)} \, \xi
    \end{equation*}
    which implies as in~(\ref{eq:RR-1}) that, for $B\in{\cal F}_t$,
    \begin{equation*}
      \lim_{\theta\rightarrow\infty}\PP(B\mid \langle X_{t+\theta,i}, 1 \rangle>0)
      =\lim_{\theta\rightarrow\infty}\PP(B\mid \langle X_{t+\theta},{\bf 1} \rangle>0)
      =\PP^*(B).
    \end{equation*}
    Therefore, Theorem~\ref{thm:cond} remains valid if the
    conditioning by the non-extinction of the whole population
    is replaced by the non-extinction of type $i$ only. This
    property relies strongly on the irreducibility of the mutation
    matrix $D$. In Section~\ref{sec:particular-case}, we will show
    that it does not always hold true when $D$ is reducible (see for
    example Theorem \ref{thm:cond-CP} or Theorem \ref{cond-CP-2}).}
  \hfill $\diamondsuit$
\end{rem}

\subsection{\mbox{Laplace functional of $\PP^*$ and Martingale Problem}}
\label{sec:immigration}

To better understand the properties of $\PP^*$, its Laplace functional
provides a very useful tool.

\begin{thm}
  \label{thm:V}
  $\PP^*$ is characterized by: $\forall f\in  C_b(\RR,\RR^k)_+$
  \begin{equation}
    \label{eq:TL-V}
    \EE^*(\exp-\langle X_t,f\rangle\mid X_0=m)
    =\frac{\langle m,V_tf\rangle}{\langle m,\xi\rangle}e^{-\mu
      t}e^{-\langle m,U_tf\rangle}
  \end{equation}
  where the semigroup $V_tf$ is the unique solution of the PDE
  \begin{equation}
    \label{eq:def-V}
    \frac{\partial V_tf}{\partial t}=\Delta V_tf+DV_tf-cU_tf\odot
    V_tf,\quad V_0f=\xi.
  \end{equation}
\end{thm}

\bigskip

\paragraph{Proof}
From Theorem~\ref{thm:cond} and (\ref{eq:LaplaceMDW}) we get
\begin{align*}
  \EE^*(e^{-\langle X_t,f\rangle}\mid X_0=m) & =\EE\left(\frac{\langle
      X_t,\xi\rangle}{\langle m, \xi\rangle}e^{-\mu t}e^{-\langle
      X_t,f\rangle}\mid X_0=m \right) \\
  & =\frac{e^{-\mu t}}{\langle
    m,\xi\rangle}\frac{\partial}{\partial\varepsilon}
  \EE\left(e^{-\langle
      X_t,f+\varepsilon\xi\rangle}\right)\bigr|_{\varepsilon=0} \\
  & =\frac{e^{-\mu t}}{\langle m,\xi\rangle}e^{-\langle m,U_tf\rangle}
  \frac{\partial}{\partial\varepsilon}\langle
  m,U_t(f+\varepsilon\xi)\rangle\bigr|_{\varepsilon=0}.
\end{align*}
Let $V_tf:=\frac{\partial}{\partial\varepsilon}
U_t(f+\varepsilon\xi)\bigr|_{\varepsilon=0}$. Then $V_tf$ is solution
of
\begin{align*}
  \frac{\partial V_tf}{\partial t} & =\frac{\partial}{\partial\varepsilon}
  \Big( \Delta U_f(f+\varepsilon\xi)+DU_t(f+\varepsilon\xi)
    -\frac{c}{2}U_t(f+\varepsilon\xi)^{\odot 2}\Big)\Bigr|_{\varepsilon=0}
  \\ & =(\Delta +D)V_tf-cU_tf\odot V_tf
\end{align*}
and $V_0f=\frac{\partial}{\partial\varepsilon}
(f+\varepsilon\xi)\bigr|_{\varepsilon=0}=\xi$.\hfill$\Box$\bigskip

Comparing with the Laplace functional of $\PP$, the multiplicative
term $\frac{\langle m,V_tf\rangle}{\langle m,\xi\rangle}e^{-\mu t}$
appears in the Laplace functional of $\PP^*$. 
In particular, the multitype Feller diffusion $x_t$ is characterized under $\PP^*$ by
\begin{equation}
  \label{eq:TL-v}
  \EE^*(\exp-( x_t,\lambda)\mid x_0=x)
  =\frac{( x,v_t^\lambda)}{( x,\xi)}e^{-\mu
    t}e^{-( x,u_t^\lambda)}, \quad  \lambda\in\RR^k_+
\end{equation}
where $v_t^\lambda:=V_t\lambda$ satisfies the differential system
\begin{equation}
  \label{eq:def-v}
  \frac{dv_t^\lambda}{dt}=Dv_t^\lambda-cu_t^\lambda\odot
  v_t^\lambda,\quad v_0^\lambda=\xi.
\end{equation}

The following theorem gives the martingale problem satisfied by
the conditioned MDW process. This formulation also allows one to
interpret $\PP^*$ as an unconditioned MDW process with immigration
(see Remark~\ref{rem:immigr} below). The term with Laplace
functional $\frac{\langle m,V_tf\rangle}{\langle m,\xi\rangle}e^{-\mu
  t}$ that we just mentioned is another way to interpret this
immigration.

\begin{thm}
  \label{thm:immigration}
  $\PP^*$ is the unique solution of the following martingale problem:
  for all $f\in C_b^2(\RR,\RR^k)_+$,
  \begin{multline} \label{eq:mart-pb}
     \exp(-\langle X_t,f\rangle)-\exp(-\langle m, f\rangle) \\
    + \int_0^t\Big( \langle X_s,(\Delta + D )f\rangle
      + c \frac{\langle X_s,f\odot\xi\rangle}{\langle X_s,\xi\rangle}
      -\frac{c}{2}\langle X_s,f^{\odot 2}\rangle\Big) \exp(-\langle
    X_s,f\rangle) \, ds
  \end{multline}
  is a $\PP^*$-local martingale.
\end{thm}

\paragraph{Proof}
According to \cite{GL-M90AAP} (see also \cite{EK-RC91SPA} for the
monotype case), $\PP$ is the unique solution of the following
martingale problem~:
for any function $F:M(\RR)^k \rightarrow \RR$ of the form
$\varphi(\langle\cdot,f\rangle)$ with $\varphi\in C^2(\RR,\RR)$ and
$f\in C^2_b(\RR,\RR^k)_+$,
\begin{equation}
\label{eq:mart-pb-P}
  F(X_t)-F(X_0)-\int_0^t{\cal A}F(X_s) \, ds \quad \textrm{is a
    $\PP$-local martingale}.
\end{equation}
 Here the infinitesimal generator ${\cal A}$ is  given by
\begin{align*}
  {\cal A}F(m) & =\langle m,(\Delta +D)\frac{\partial F}{\partial
    m}\rangle+\frac{c}{2}\langle m,\partial^2 F/\partial m^2\rangle \\
  & = \sum_{i=1}^k\langle m_i,\Delta \frac{\partial F}{\partial
    m_i}+\sum_{j=1}^kd_{ij}\frac{\partial F}{\partial m_j}\rangle
  +\frac{c}{2}\sum_{i=1}^k\langle m_i,\frac{\partial^2F}{\partial m_i^2}\rangle .
\end{align*}
where we use the notation $\partial F/\partial m=(\partial F/\partial
m_i)_{1\leq i\leq k}$ and $\partial^2 F/\partial m^2=(\partial^2
F/\partial m_i^2)_{1\leq i\leq k}$ with
\begin{equation*}
  \frac{\partial F}{\partial m_i}(x):=\lim_{\varepsilon\rightarrow 0}
  \frac{1}{\varepsilon}\big(
  F(m_1,\ldots,m_i+\varepsilon\delta_x,\ldots,m_k)-F(m) \big), x \in
  \RR.
\end{equation*}
Applying this to the time-dependent function
$$
F(s,m):=\langle
m,\xi\rangle e^{- \mu s}e^{-\langle m,f\rangle} \textrm{ with } f\in
C^2_b(\RR,\RR^k)_+
$$
for which
\begin{align*}
  \frac{\partial F(s,m)}{\partial m}(x) &
  =-f(x)F(s,m)+\xi e^{-\mu s -\langle m,f\rangle} \\
  \textrm{ and  } \quad \frac{\partial^2 F(s,m)}{\partial m^2} (x) &
  =f^{\odot 2}(x)F(s,m)-2 f(x)\odot\xi e^{-\mu s -\langle m,f\rangle}, 
\end{align*}
one gets
\begin{multline*}
  \frac{\partial F}{\partial s}(s,m)+{\cal A}(F(s,\cdot))(m) \\
  =-\langle m,(\Delta+D)f\rangle F+\frac{c}{2}\langle m,f^{\odot
    2}\rangle F-c\frac{\langle m,f\odot\xi\rangle}{\langle m,\xi\rangle}F.
\end{multline*}
Therefore,
\begin{multline*}
  \langle X_t,\xi\rangle e^{-\mu t-\langle X_t,f\rangle}-\langle
  m,\xi\rangle e^{-\langle m, f\rangle} \\ + \int_0^t\langle
  X_s,\xi\rangle e^{-\mu s}\Big( \langle X_s,(\Delta +D)f\rangle
  +c\frac{\langle X_s,f\odot\xi\rangle}{\langle X_s,\xi\rangle}
  -\frac{c}{2}\langle X_s,f^{\odot 2}\rangle\Big) e^{-\langle X_s,f\rangle}ds
\end{multline*}
is a $\PP$-local martingale, which implies that~(\ref{eq:mart-pb}) is
a $\PP^*$-local martingale.\\
The uniqueness of the solution $\PP^*$ to the martingale
problem~(\ref{eq:mart-pb}) comes from the uniqueness of the solution
of the martingale problem~(\ref{eq:mart-pb-P}).\hfill$\Box$\bigskip

\begin{rem}
  \label{rem:immigr}
  \textup{Due to the form of the martingale problem
    (\ref{eq:mart-pb}), the probability measure $\PP^*$ can be
    interpreted as the law of a MDW process with interactive
    immigration whose rate at time $s$, if conditioned by $X_s$, is a
    random measure with Laplace functional $ \exp-c\frac{\langle
      X_s,f\odot\xi\rangle}{\langle X_s,\xi\rangle}$. Monotype DW
    processes with deterministic immigration rate were introduced by
    Dawson in~\cite{Da78Bolyai}. The first interpretation of
    conditioned branching processes as branching processes with
    immigration goes back to Kawazu and Watanabe in~\cite{KaWa71TPA},
    Example~2.1. See also \cite{RR89CRAS} and \cite{EvPe} for further
    properties in the monotype case.}\hfill $\diamondsuit$
\end{rem}
\bigskip

\section{Long time behavior of conditioned multitype Feller diffusions}
\label{sec:long-time}

We are now interested in the long time behavior of the MDW process
conditioned on non-extinction in the remote future. Unfortunately,
because of the Laplacian term in~(\ref{eq:def-V}), there is no hope to
obtain a limit of $X_t$ under $\PP^*$ at the level of measure
(however, see~\cite{EvPe} for the long time behavior of critical
monotype conditioned Dawson-Watanabe processes with \emph{ergodic}
spatial motion). Therefore, we will restrict our attention to the
$\RR^k$-valued multitype Feller diffusion $x_t$.  As a first step in
our study, we begin this section with the monotype case.

\subsection{Monotype case}
\label{sec:long-time-monotype}

In this subsection, we first give asymptotic behavior of $x_t$
under $\PP^*$ (Proposition~\ref{thm:t-infty-monotype}). This result is
already known, but we give a proof that will be useful in the
following section. We also give a new result about the exchange of
limits (Proposition~\ref{prop:interv-CP-monotype}).

Let us first introduce some notation for the monotype case.\\
The matrix $D$ is reduced to its eigenvalue $\mu$, the vector $\xi$ is
equal to the number 1.  Since we only consider the critical and
subcritical cases, one has $\mu \le 0$. The law $\PP^*$ of the MDW
process conditioned on non-extinction in the remote future is locally
absolutely continuous with respect to $\PP$ (monotype version of
Theorem \ref{thm:cond}, already proved in \cite{RR89CRAS}, Proposition
1).  More precisely
\begin{equation}
  \label{eq:cond-monotype}
  \PP^*\bigr|_{{\cal F}_t}=\frac{\langle X_t,1\rangle}{\langle
    m, 1\rangle}e^{- \mu  t}\, \PP\bigr|_{{\cal F}_t}.
\end{equation}
Furthermore the Laplace functional of $\PP^*$ satisfies
(see~\cite{RR89CRAS}, Theorem 3):
\begin{equation}
  \label{eq:TL-V-monotype}
  \EE^*(\exp-\langle X_t,f\rangle\mid X_0=m)
  =\frac{\langle m,V_tf\rangle}{\langle m, 1 \rangle} e^{- \mu
    t}e^{-\langle m,U_tf\rangle}
  =\frac{\langle m, \tilde V_tf\rangle}{\langle m, 1 \rangle}
  e^{-\langle m,U_tf\rangle}
\end{equation}
where
\begin{equation}
  \label{eq:def-V-monotype}
  \frac{\partial \tilde V_tf}{\partial t}=\Delta \tilde V_tf -c\, U_tf 
  \tilde V_tf,\quad \tilde V_0f=1.
\end{equation}

The total mass process $x_t=\langle X_t,1 \rangle$ is a (sub)critical
Feller branching diffusion under $\PP$. By~(\ref{eq:TL-V-monotype}) its
Laplace transform under $\PP^*$ is
\begin{equation}
  \label{eq:TL-v-monotype}
  \EE^*(\exp- \lambda x_t\mid x_0=x)
  =\tilde v_t^\lambda \, e^{- x u_t^\lambda}, \qquad \lambda\in\RR_+,
\end{equation}
with
\begin{equation*}
  \label{eq:def-v-monotype}
  \frac{d \tilde v_t^\lambda}{dt}=-cu_t^\lambda  \tilde
  v_t^\lambda,\quad \tilde v_0^\lambda=1.
\end{equation*}
Recall that the cumulant $u_t^\lambda$ satisfies
\begin{equation}
  \label{eq:def-u-monotype}
  \frac{d u_t^\lambda}{dt}=\mu \ u_t^\lambda -\frac{c}{2} \
  (u_t^\lambda)^2 ,\quad u_0^\lambda=\lambda.
\end{equation}
This yields in the subcritical case the explicit formulas
\begin{equation}
  \label{u-monotype}
  u_t^\lambda=\frac{\lambda  \, e^{\mu t}}{1 +\frac{c}{2|\mu|}\lambda
    (1- e^{\mu t})} , \quad \lambda\geq 0,
\end{equation}
\begin{equation}
  \label{eq:def-v-int-monotype}
 \textrm{ and } \quad  \tilde v_t^\lambda=\exp\left(-c\int_0^t u_s^\lambda ds\right)
   = \frac{1}{\big(1 + \frac{c}{2|\mu|}\lambda (1- e^{\mu t})\big)^2}.
\end{equation}
In the critical case ($\mu=0$) one obtains (see~\cite{LN68TPA} Equation~(2.14))
\begin{equation}
  \label{u-v-monotype-critique}
    u_t^\lambda=\frac{\lambda }{1+\frac{c}{2}\lambda t} \qquad \textrm{ and }
    \qquad v_t^\lambda= \tilde v_t^\lambda = \frac{1}{\big(1 +
      \frac{c}{2}\lambda t\big)^2}.
\end{equation}

We are now ready to state the following asymptotic result. 
\begin{prop}
  \label{thm:t-infty-monotype} \hspace{10cm}
  \begin{description}
  \item[\textmd{(a)}] In the critical case ($\mu=0$), the process
    $x_t$ explodes in $\PP^*$-probability when $t\rightarrow \infty$,
    i.e. for any $M> 0$,
    $$ 
    \lim_{t\rightarrow +\infty}\PP^*(x_t\leq M) = 0.
    $$
  \item[\textmd{(b)}] In the subcritical case ($\mu<0$),
    $$
    \lim_{t\rightarrow +\infty}  \PP^*(x_t \in \cdot ) \overset{(d)}{=} \Gamma(2, \frac{2|\mu|}{c})
    $$
    where this notation means that $x_t$ converges in
    $\PP^*$-distribution to a Gamma distribution with parameters $2$ and $2|\mu|/c$.
  \end{description}
\end{prop}

One can find in~\cite{Lambert06} Theorem~4.2 a proof of this theorem
for a more general model, based on a pathwise approach.  We propose
here a different proof, based on the behavior of the cumulant
semigroup and moment properties, which will be useful in the sequel.

\paragraph{Proof}
For $\mu=0$, by (\ref{u-v-monotype-critique}), $u_t^\lambda\rightarrow
0$ and $v_t^\lambda\rightarrow 0$ when $t\rightarrow\infty$ for any
$\lambda\not =0$. This implies by (\ref{eq:TL-v-monotype}) the
asymptotic explosion of $x_t$ in $\PP^*$-probability.  Actually, the
rate of explosion is also known: in \cite{EvPe} Lemma 2.1, the authors
have proved that $\frac{x_t}{t}$ converges in distribution as $t
\rightarrow \infty$ to a Gamma-distribution. This can also be deduced
from~(\ref{u-v-monotype-critique}), since $u^{\lambda/t}_t$ and
$v^{\lambda/t}_t$ converge
to 0 and $1/(1+c\lambda/2)^2$ respectively, as $t\rightarrow\infty$.\\

For $\mu<0$, by (\ref{eq:TL-v-monotype}), (\ref{u-monotype}) and
(\ref{eq:def-v-int-monotype}), the process $x_t$ has the
same law as the sum of two independent random variables, the first one
with distribution $\Gamma(2,\frac{2|\mu|}{c(1- e^{\mu t})})$ and the
second one vanishing for $t\rightarrow\infty$.
The conclusion is now clear.\hfill$\Box$\bigskip
    
\begin{rem}
  \textup{The presence of a Gamma-distribution in the above
    Proposition is not surprising.
    \begin{itemize}
    \item As we already mentioned it appears in the critical case as
      the limit law of $x_t/t$ \cite{EvPe}.
    \item It also goes along with the fact that these distributions
      are the equilibrium distributions for subcritical Feller
      branching diffusions with constant immigration.  (See
      \cite{AN72}, and Lemma 6.2.2 in \cite{DGWEJP04}).  \emph{We are
        grateful to A. Wakolbinger for proposing this interpretation.}
    \item Another interpretation is given in~\cite{Lambert06}. The
      Yaglom distribution of the process $x_t$, defined as the limit
      law as $t\rightarrow\infty$ of $x_t$ conditioned on $x_t>0$, is
      the exponential distribution with parameter $2|\mu|/c$ (see
      Proposition~\ref{prop:interv-CP-monotype} below, with
      $\theta=0$). The Gamma distribution appears as the size-biased
      distribution of the Yaglom limit ($\PP^*(x_\infty\in
      dr)=r\PP(Y\in dr)/\EE(Y)$, where $Y\sim {\cal
        E}\mbox{\textit{xp}}(\frac{2|\mu|}{c})$), which is actually
      a general fact (\cite[Th.4.2(ii)(b)]{Lambert06}).  \hfill
      $\diamondsuit$
    \end{itemize}}
\end{rem}
\bigskip
    
We have just proved that, for $\mu <0$, the law of $x_t$ conditioned
on $x_{t+\theta}>0$ converges to a Gamma distribution when taking
first the limit $\theta\rightarrow\infty$ and next the limit
$t\rightarrow\infty$. It is then natural to ask whether the order of
the two limits can be exchanged: what happens if one first fix $\theta$
and let $t$ tend to infinity, and then let $\theta$ increase?  We
obtain the following answer.

\begin{prop}
  \label{prop:interv-CP-monotype}
  When $\mu <0$, conditionally on $x_{t+\theta}>0$, $x_t$ converges in
  distribution when $t\rightarrow\infty$ to the sum of two independent
  exponential r.v.\ with respective parameters $\frac{2|\mu|}{c}$ and
  $\frac{2|\mu|}{c}(1-e^{ \mu \theta})$.
  
  Therefore, one can interchange both limits in time:
  \begin{equation*}
    \lim_{\theta\rightarrow\infty}\lim_{t\rightarrow\infty}
    \PP(x_{t}\in\cdot \mid x_{t+\theta}>0)
    \overset{(d)}{=} \lim_{t\rightarrow\infty}\lim_{\theta\rightarrow\infty}
    \PP(x_{t}\in\cdot \mid x_{t+\theta}>0) 
    \overset{(d)}{=}\Gamma(2,\frac{2|\mu|}{c}).
  \end{equation*}
\end{prop}

\paragraph{Proof} 
First, observe that, by~(\ref{u-monotype}),
$$
\lim_{\bar\lambda\rightarrow\infty}u_t^{\bar\lambda}=\frac{2|\mu|}{c}
\frac{e^{\mu t}}{1-e^{\mu t}} \textrm{ and }
\lim_{t\rightarrow\infty}\frac{ u_t^{\lambda}}{e^{\mu t}}
=\frac{\lambda}{1-\frac{c}{2\mu} \lambda }.
$$
As in~(\ref{eq:RR-1}), it holds 
\begin{align*}
  \EE(e^{-\lambda x_t}\mid x_{t+\theta}>0) &
  =\frac{\EE\big(e^{-\lambda x_t}(1-\PP( x_{t+\theta}=0\mid{\cal F}_t))\big)}
  {1-\PP(x_{t+\theta}=0)}  \\
  & =\frac{\EE\big(e^{-\lambda x_t}(1-
    e^{-x_t\lim_{\bar\lambda\rightarrow\infty}u_\theta^{\bar\lambda}})
    \big)}
  {1- e^{- x \lim_{\bar\lambda\rightarrow\infty}u_{t+\theta}^{\bar\lambda}} } \\
  & =\frac{e^{- x u_t^{\lambda}}- e^{-x u_t^{\lambda
        +\lim_{\bar\lambda\rightarrow\infty}u_\theta^{\bar\lambda}}} }
  {1- e^{-x \lim_{\bar\lambda\rightarrow\infty}u_{t+\theta}^{\bar\lambda}} } \\
  & =\frac{e^{- x u_t^{\lambda}}- e^{-x u_t^{\lambda
        +\lim_{\bar\lambda\rightarrow\infty}u_\theta^{\bar\lambda}}} }
  {1- \exp(- x \frac{2|\mu|}{c} \frac{e^{\mu (t+\theta) }}{1- e^{\mu (t+\theta)}}) }   
\end{align*}
Thus 
\begin{align*}
  \lim_{t \rightarrow\infty}\EE(e^{-\lambda x_t}\mid x_{t+\theta}>0) &
  = \frac{c}{2|\mu|} e^{-\mu \theta }
  \lim_{t \rightarrow\infty} 
  e^{ |\mu| t} \big( u_t^{\lambda
    +\lim_{\bar\lambda\rightarrow\infty}u_\theta^{\bar\lambda}} -
  u_t^{\lambda}\big) \\
  & =  \frac{c}{2|\mu|} e^{-\mu \theta }
  \Big(
  \frac{\lambda +\lim_{\bar\lambda\rightarrow\infty}u_\theta^{\bar\lambda}}
  {1-\frac{c}{2\mu} (\lambda
    +\lim_{\bar\lambda\rightarrow\infty}u_\theta^{\bar\lambda})}-
  \frac{\lambda}{1-\frac{c}{2\mu} \lambda }
  \Big) \\
  & = \frac{1}{ 1 + \frac{c}{2|\mu|} \lambda }\cdot  \frac{1}{
    1+\frac{c}{2|\mu|}(1- e^{\mu \theta })\lambda }
\end{align*}
where the first (resp.\ the second) factor is equal to the Laplace
transform of an exponential r.v.\ with parameter $2 |\mu|/c$ (resp.\
with parameter $\frac{2|\mu|}{c}(1-e^{\mu \theta })$).  This means
that
\begin{equation*}
  \lim_{t \rightarrow\infty}  \PP(x_{t}\in\cdot \mid x_{t+\theta}>0)
  \overset{(d)}{=} {\cal E}\mbox{\textit{xp}}(\frac{2|\mu|}{c})\otimes{\cal E}
  \mbox{\textit{xp}}(\frac{2|\mu|}{c}(1-e^{\mu \theta })).
\end{equation*}
It is now clear that 
\begin{equation*}
  \lim_{\theta \rightarrow\infty}  \lim_{t \rightarrow\infty}
  \PP(x_{t}\in\cdot \mid x_{t+\theta}>0)
  \overset{(d)}{=}{\cal E}\mbox{\textit{xp}}(\frac{2|\mu|}{c})\otimes{\cal E}
  \mbox{\textit{xp}}(\frac{2|\mu|}{c}) = \Gamma(2,\frac{2|\mu|}{c}).
\end{equation*} 
Thus, the limits in time interchange.\hfill$\Box$\bigskip

\begin{rem}
  \label{rem:d=0-monotype}
  \textup{The previous computation is also possible in the critical
    case and gives a similar interchangeability result. More
    precisely, for any $\theta>0$, $x_t$ explodes conditionally on
    $x_{t+\theta}>0$ in $\PP$-probability when
    $t\rightarrow+\infty$. In particular, for any $M>0$,
    \begin{equation*}
      \lim_{\theta\rightarrow\infty}\lim_{t\rightarrow\infty}
      \PP(x_{t}\leq M\mid x_{t+\theta}>0)
      = \lim_{t\rightarrow\infty}\lim_{\theta\rightarrow\infty}
      \PP(x_{t}\leq M\mid x_{t+\theta}>0)=0.
    \end{equation*} }
  \hfill $\diamondsuit$
\end{rem}

\begin{rem}
  \label{rem:beta-stable}
  \textup{One can develop the same ideas as before when the branching
    mechanism with finite variance $c$ is replaced by a $\beta$-stable
    branching mechanism, $0<\beta<1$, with infinite variance
    (see~\textup{\cite{Da78Bolyai}} Section 5 for a precise
    definition).  In this case,
    equation~\textup{(\ref{eq:def-u-monotype})} has to be replaced by
    $ \frac{du_t^\lambda}{dt}=\mu \, u_t^\lambda -
    c(u_t^\lambda)^{1+\beta}$ which implies that
    \begin{equation*}
      u_t^\lambda=\frac{\lambda \,
        e^{\mu t}}{\left( 1+\frac{c\lambda^\beta}{|\mu|}(1- e^{\beta
            \mu t}) \right)^{1/\beta}}.
    \end{equation*}
    Therefore, with a similar calculation as above, one can easily
    compute the Laplace transform of the limit conditional law of
    $x_t$ when $t\rightarrow\infty$ and prove the exchangeability
    of limits:
    \begin{equation*}
      \lim_{\theta\rightarrow\infty}\lim_{t\rightarrow\infty}
      \EE(e^{-\lambda x_t}\mid x_{t+\theta}>0)
      =\lim_{t\rightarrow\infty}\lim_{\theta\rightarrow\infty}
      \EE(e^{-\lambda x_t}\mid x_{t+\theta}>0)
      =\frac{1}{\left(1+\frac{c}{|\mu|}\lambda^\beta\right)^{1+1/\beta}}.
    \end{equation*}
    As before, this distribution can be interpreted as the
    size-biased Yaglom distribution corresponding to the stable
    branching mechanism.  This conditional limit law for the
    subcritical branching process has been obtained in~\cite{Li00}
    Theorem 4.2. We also refer to~\cite{Lambert06} Theorem 5.2 for a
    study of the critical stable branching process. } \hfill
  $\diamondsuit$
\end{rem}

\subsection{Multitype irreducible case}
\label{sec:long-time-multitye}

We now present the multitype generalization of Proposition
\ref{thm:t-infty-monotype} on the asymptotic behavior of the 
conditioned multitype Feller diffusion with irreducible mutation matrix $D$.

\begin{thm}
  \label{thm:t-infty}
  \begin{description}
  \item[\textmd{(a)}] In the critical case ($\mu=0$), when the
    mutation matrix $D$ is irreducible, $x_t$
    explodes in $\PP^*$-probability when $t\rightarrow \infty$, i.e.
    $$
    \forall i\in\{1,\ldots,k\} ,\,  \forall M>0, \quad
    \lim_{t\rightarrow+\infty}\PP^*(x_{t,i}\leq M)=0.
    $$
  \item[\textmd{(b)}] In the subcritical case ($\mu<0$) when $D$ is
    irreducible, the law of $x_t$ converges in
    distribution under $\PP^*$ when $t\rightarrow \infty$ to a
    non-trivial limit which does not depend on the initial condition
    $x$.
  \end{description}
\end{thm}

\paragraph{Proof}

One obtains from~(\ref{eq:def-v}) that
\begin{equation*}
  Dv_t^\lambda-c\sup_i(u_{t,i}^\lambda) \,  v_t^\lambda
  \leq\frac{dv_t^\lambda}{dt}
  \leq Dv_t^\lambda-c\inf_i(u_{t,i}^\lambda) \, v_t^\lambda.
\end{equation*}
Then, by Lemma~\ref{lem:diff-ineq},
\begin{equation}
  \label{eq:ineq-v}
  \exp\Big(\mu t-c\int_0^t\sup_iu_{s,i}^\lambda\: ds\Big)\xi\leq
  v_t^\lambda\leq\exp\Big(\mu t-c\int_0^t\inf_iu_{s,i}^\lambda\: ds\Big)\xi.
\end{equation}
Therefore, in the critical case, $v_t^\lambda$ vanishes for $t$ large
if $\lambda>0$, due to the divergence of $\int_0^\infty
\inf_iu_{s,i}^\lambda ds$, which is itself a consequence of
Lemma~\ref{lem:pte-u} (iii).  If $\lambda_i=0$ for some type $i$,
by Lemma~\ref{lem:pte-u} (i) and the semigroup property of $t \mapsto
u_t$, we can use once again Lemma~\ref{lem:pte-u} (iii) starting from
a positive time, to prove that $\lim_{t\rightarrow\infty} v_t^\lambda=
0$.  Then, the explosion of $x_t$ in $\PP^*$-probability follows
directly from~(\ref{eq:TL-v}) and from the fact that
$\lim_{t\rightarrow\infty}u_t^\lambda= 0$.

To prove~(b), we study the convergence of ${\tilde
  v}_t^\lambda:=e^{-\mu t}v_t^\lambda$ when $t\rightarrow\infty$.
By~(\ref{eq:ineq-v}) and Lemma~\ref{lem:pte-u} (ii) we know that
$t\mapsto {\tilde v}_t^\lambda$ is bounded and bounded away from $0$.
Fix $\varepsilon\in(0,1)$ and $t_0$ such that
$\int_{t_0}^\infty\sup_iu_{t,i}^\lambda \:dt<\varepsilon$. Then, for
any $t\geq 0$,
\begin{equation}
  \label{eq:indep-CI}
  e^{-c\varepsilon}e^{(D-\mu I)t}{\tilde v}_{t_0}^\lambda\leq
  {\tilde v}_{t_0+t}^\lambda\leq e^{(D-\mu I)t}{\tilde v}_{t_0}^\lambda
\end{equation}
and so, for any $s,t\geq 0$,
\begin{equation*}
  |{\tilde v}_{t_0+t+s}^\lambda-{\tilde v}_{t_0+t}^\lambda|\leq
  \sup_{\delta\in\{-1,1\}}\left|\left(e^{c\delta\varepsilon}I-e^{(D-\mu
        I)s}\right)e^{(D-\mu I)t}{\tilde v}_{t_0}^\lambda\right|.
\end{equation*}
By Perron-Frobenius' theorem, $\lim_{t\rightarrow\infty}e^{(D-\mu
  I)t}= P:=(\xi_i\eta_j)_{i,j}$ and thus, when $t\rightarrow\infty$,
\begin{equation*}
  |{\tilde v}_{t_0+t+s}^\lambda-{\tilde v}_{t_0+t}^\lambda|\leq (e^{c\varepsilon}-1)
  ({\tilde v}_{t_0}^\lambda,\eta)\, \xi+ |{\tilde v}_{t_0}^\lambda|o(1)
\end{equation*}
where the negligible term $o(1)$ does not depend on $t_0,\ s,\
\varepsilon$ and $\lambda$, since $\varepsilon<1$ and $\exp((D-\mu
I)s)$ is a bounded function of $s$.  Therefore, $({\tilde
  v}_t^\lambda)_{t\geq 0}$ satisfies the Cauchy criterion and
converges to a finite positive limit ${\tilde v}^\lambda_\infty$ when
$t\rightarrow\infty$.

We just proved the convergence of the Laplace functional~(\ref{eq:TL-v}) of $x_t$
under $\PP^*$ when $t\rightarrow\infty$.  In order to obtain the
convergence in law of $x_t$, we have to check the continuity of the
limit for $\lambda=0$, but this is an immediate consequence of
$\lim_{\lambda\rightarrow 0}\lim_{t\rightarrow\infty}{\tilde
  v}_t^\lambda=\xi$.

Finally, letting $t$ go to infinity in~(\ref{eq:indep-CI}), we get
$$
|{\tilde v}^\lambda_\infty-P{\tilde v}^\lambda_{t_0}|\leq
(1-e^{-c\varepsilon})P{\tilde v}^\lambda_{t_0},
$$
where $P{\tilde v}^\lambda_{t_0}=({\tilde v}^\lambda_{t_0},\eta)\xi$.
It follows that ${\tilde v}^\lambda_\infty$ is proportional to $\xi$,
as limit of quantities proportional to $\xi$. Therefore $(x,{\tilde
  v}^\lambda_\infty)/(x,\xi)= ({\tilde v}^\lambda_\infty,{\bf 1})$ is
independent of $x$ and the limit law of $x_t$ too.\hfill$\Box$\bigskip

We can also generalize the exchange of limits of
Proposition~\ref{prop:interv-CP-monotype} to the multitype irreducible
case.

\begin{thm}
  \label{thm:interv-multi}
  In the subcritical case, conditionally on $(x_{t+\theta},{\bf
    1})>0$, $x_t$ converges in distribution when $t\rightarrow+\infty$
  to a non-trivial limit which depends only on $\theta$. Furthermore,
  one can interchange both limits in $t$ and $\theta$ :
  \begin{equation*}
    \lim_{\theta\rightarrow\infty}\lim_{t\rightarrow\infty}
    \PP(x_{t}\in\cdot \mid (x_{t+\theta},{\bf 1})>0)
    \overset{(d)}{=} \lim_{t\rightarrow\infty}\lim_{\theta\rightarrow\infty}
    \PP(x_{t}\in\cdot \mid (x_{t+\theta},{\bf 1})>0).
  \end{equation*}
\end{thm}

\paragraph{Proof}
Following a similar computation as in the proof of
Proposition~\ref{prop:interv-CP-monotype}, 
\begin{align*}
  \lim_{t\rightarrow\infty}\EE(e^{-(x_t,\lambda)}\mid
  (x_{t+\theta},{\bf 1})>0) &
  =\lim_{t\rightarrow\infty}\frac{\exp(-(x,u^\lambda_t))
    -\exp(-(x,u_t^{\lambda+\lim_{\bar\lambda\hookrightarrow\infty}u^{\bar\lambda}_\theta}))}
  {1-\exp(-(x,\lim_{\bar\lambda\hookrightarrow\infty}u^{\bar\lambda}_{t+\theta}))} \\
  & =\lim_{t\rightarrow\infty}
  \frac{(x,u^{\lambda+\lim_{\bar\lambda\hookrightarrow\infty}u^{\bar\lambda}_\theta}_t
    -u^\lambda_t)}{(x,\lim_{\bar\lambda\hookrightarrow\infty}u^{\bar\lambda}_{t+\theta})} .
\end{align*}
Since 
$$
D u_t^\lambda -\frac{c}{2}(\sup_iu_{t,i}^\lambda) \,  u_t^\lambda \leq
\frac{du_t^\lambda}{dt} \leq 
D u_t^\lambda  -\frac{c}{2}(\inf_iu_{t,i}^\lambda) \, u_t^\lambda ,
$$
 one gets~:
\begin{equation*}
  \exp\Big(-\frac{c}{2}\int_0^t\sup_iu_{s,i}^\lambda\: ds\Big)e^{Dt}\lambda\leq
  u_t^\lambda\leq\exp\Big(-\frac{c}{2}\int_0^t\inf_iu_{s,i}^\lambda\:
  ds\Big)e^{Dt}\lambda.
\end{equation*}
This inequality, similar to~(\ref{eq:ineq-v}), can be used exactly as
in the proof of Theorem~\ref{thm:t-infty} (b) to prove that, for any
$\varepsilon$, there exists $t_0$ large enough such that
\begin{equation}
  \label{eq:indep-CI-u-infty}
  e^{-c\varepsilon/2}e^{(D-\mu I)t}e^{-\mu t_0}u_{t_0}^\lambda\leq
  e^{-\mu(t+t_0)}u_{t_0+t}^\lambda\leq e^{(D-\mu I)t}e^{-\mu t_0}u_{t_0}^\lambda .
\end{equation}
and to deduce from~(\ref{eq:indep-CI-u-infty}) that ${\tilde
  u^\lambda_t}:= e^{-\mu t}u^\lambda_t$ converges as $t\rightarrow
\infty$ to a non-zero limit ${\tilde u^\lambda_\infty}$ proportional
to the vector $\xi$.

Moreover, because of~(\ref{eq:TL-x}), $u^\lambda_t$ is increasing with
respect to each coordinate of $\lambda$. Therefore, it is elementary
to check that $t\mapsto\lim_{\bar\lambda\hookrightarrow\infty}
u^{\bar\lambda}_t=\sup_n u^{n\mathbf{1}}_t$ is also solution of the
non-linear differential system~(\ref{eq:ODE-u}), but only defined on
$(0,\infty)$ (recall that, by Lemma~\ref{lem:pte-u} (ii),
$\lim_{\bar\lambda\hookrightarrow\infty}u^{\bar\lambda}_t<\infty$ for
any $t>0$). Indeed, assume that $b_t=\sup_n a^n_t$ where
$\dot{a}^n_t=F(a^n_t)$ for a locally Lipschitz function $F$. Fix $t$
such that $b_t<+\infty$ and a small $\eta>0$, and let
$\underline{F}:=\inf_{|x-b_t|\leq\eta}F(x)$ and
$\bar{F}:=\sup_{|x-b_t|\leq\eta}F(x)$. There exists $n_0$ such that,
for $n\geq n_0$, $|a^n_t-b_t|\leq \eta/2$. Moreover, for any $s$ in a
neighborhood of $t$, $|a^n_s-b_t|\leq\eta$, where the neighborhood
depends on $\bar{F}$ and $\underline{F}$, but is uniform in $n\geq
n_0$. Therefore, for sufficiently small $s$ and for $n$ sufficiently
large, $\underline{F}\leq (a^n_{t+s}-a^n_t)/s\leq \bar{F}$. Letting
$n\rightarrow\infty$, $s\rightarrow 0$ and finally $\eta\rightarrow
0$, since $\bar{F}-\underline{F}\rightarrow 0$ when $\eta\rightarrow
0$, $b_t$ is differentiable at time $t$ and $\dot{b}_t=F(b_t)$.

Therefore, the semigroup property of the flow of~(\ref{eq:ODE-u})
implies that, for any $t\geq 0$,
$$
\lim_{\bar\lambda\hookrightarrow\infty}u^{\bar\lambda}_{t+\theta}
=u_t^{\lim_{\bar\lambda\hookrightarrow\infty}u^{\bar\lambda}_\theta},
$$
so that $e^{-\mu
  t}\lim_{\bar\lambda\hookrightarrow\infty}u^{\bar\lambda}_t$ also
converges as $t\rightarrow \infty$ to a positive limit ${\tilde
  u^\infty_\infty}$, proportional to $\xi$ too.

Hence,
\begin{align}
  \lim_{t\rightarrow\infty}\EE(e^{-(x_t,\lambda)}\mid
  (x_{t+\theta},1)>0)
  & =e^{-\mu\theta} \frac{(x,{\tilde u^{\lambda+\lim_{\bar\lambda\hookrightarrow\infty}
        u^{\bar\lambda}_\theta}_\infty} -{\tilde
      u^\lambda_\infty})}{(x,{\tilde u^\infty_\infty} )} \nonumber \\
  & =e^{-\mu\theta} \frac{({\tilde u^{\lambda+\lim_{\bar\lambda\hookrightarrow\infty}
        u^{\bar\lambda}_\theta}_\infty} -{\tilde
      u^\lambda_\infty},{\bf 1})}{({\tilde u^\infty_\infty},{\bf 1})},
  \label{eq:cvgce-Lapl-trans}
\end{align}
which is independent of the initial condition $x$.

In order to prove the convergence in distribution as
$t\rightarrow\infty$ of $x_t$ conditionally on $(x_{t+\theta},{\bf
  1})>0$ to a random variable with Laplace
transform~(\ref{eq:cvgce-Lapl-trans}), it remains to prove the
continuity of this expression as a function of $\lambda$ for
$\lambda\rightarrow 0$.  To this aim and also to prove the exchangeability
of limits, we use the following Lemma, the proof of which is postponed at
the end of the subsection. This lemma gives the main reason why the
limits can be exchanged: $v^\lambda_t$ is solution of the
\emph{linearized} equation of~(\ref{eq:ODE-u}), and therefore, the
gradient of $u^\lambda_t$ with respect to $\lambda$ is solution of the
\emph{same} system of ODEs as $v^\lambda_t$. The function
$v^\lambda_t$ was involved in the computation of
$\lim_t\lim_\theta\PP(x_t\in\cdot\mid x_{t+\theta}>0)$, whereas the
gradient of $u^\lambda_t$ will be involved in the computation of
$\lim_\theta\lim_t\PP(x_t\in\cdot\mid x_{t+\theta}>0)$.

\begin{lem}
  \label{lem:u-infty-diff}
  The function $\lambda \mapsto u^{\lambda}_t$ is differentiable, and
  its derivative in the direction $\eta$, denoted by
  $\triangledown_\eta u^\lambda_t$, is solution of the same
  differential system \textup{(\ref{eq:def-v})} as $v_t^\lambda$
  except for the initial condition given by $\triangledown_\eta
  u^\lambda_0=\eta$.  Furthermore, $\lambda \mapsto {\tilde
    u}_\infty^\lambda$ is differentiable too and its derivative in the
  direction $\eta$, denoted by $\triangledown_\eta {\tilde
    u}^\lambda_\infty$ satisfies
  $$
  \triangledown_\eta {\tilde u}^\lambda_\infty=\lim_{t\rightarrow\infty} e^{-\mu
    t}\triangledown_\eta u^\lambda_t.
  $$
\end{lem}

Since $\triangledown_\eta u^\lambda_t$ satisfies the same differential
equation as $v^\lambda_t$, in particular, $\triangledown_\xi
u^\lambda_t= v^\lambda_t$ and, with the notations of the proof of
Theorem~\ref{thm:t-infty}, $\triangledown_\xi {\tilde
  u}^\lambda_\infty={\tilde v}^\lambda_\infty$.

It also follows from the above lemma that ${\tilde u}_\infty^\lambda$
is continuous as a function of $\lambda$. As a result,
\begin{align*}
  \lim_{\lambda\rightarrow 0}\lim_{t\rightarrow\infty}\EE(e^{-(x_t,\lambda)}\mid
  (x_{t+\theta},{\bf 1})>0) & =\lim_{\lambda \rightarrow
    0}e^{-\mu\theta}\frac{({\tilde u^{\lambda+\lim_{\bar\lambda\hookrightarrow\infty}
        u^{\bar\lambda}_\theta}_\infty} -{\tilde u^\lambda_\infty},{\bf
      1})}{({\tilde u^\infty_\infty},{\bf 1})} \\
  & =  e^{-\mu\theta} \frac{({\tilde u^{\lim_{\bar\lambda\hookrightarrow\infty}
        u^{\bar\lambda}_\theta}_\infty} ,{\bf 1})}{({\tilde u^\infty_\infty},{\bf 1})} 
  =1 
\end{align*}
since
$$
{\tilde u}^{\lim_{\bar\lambda\hookrightarrow\infty}u^{\bar\lambda}_\theta}_\infty 
=\lim_{t\rightarrow\infty}e^{-\mu t}
u^{\lim_{\bar\lambda\hookrightarrow\infty}u^{\bar\lambda}_\theta}_t 
=\lim_{t\rightarrow\infty}e^{-\mu
  t}\lim_{\bar\lambda\hookrightarrow\infty}u^{\bar\lambda}_{t+\theta}
=e^{\mu\theta}{\tilde u}^\infty_\infty.  
$$

Finally, let us check that the limits in $t$ and $\theta$ can be
exchanged.  Since
$\lim_{\bar\lambda\hookrightarrow\infty}u^{\bar\lambda}_\theta\sim
e^{\mu\theta}{\tilde u}^\infty_\infty=e^{\mu\theta}(\tilde{u}^\infty_\infty,\mathbf{1})\xi$ when $\theta\rightarrow\infty$,
it follows from Lemma~\ref{lem:u-infty-diff} that
\begin{equation*}
  \lim_{\theta\rightarrow\infty}e^{-\mu\theta}\frac{({\tilde
      u^{\lambda+\lim_{\bar\lambda\hookrightarrow\infty}
        u^{\bar\lambda}_\theta}_\infty} -{\tilde u^\lambda_\infty},{\bf
      1})}{({\tilde u^\infty_\infty},{\bf 1})}
  =(\triangledown_\xi {\tilde u}^\lambda_\infty ,{\bf 1})
  =({\tilde v}^\lambda_\infty,{\bf 1})
  =\lim_{t\rightarrow\infty}\EE^*(e^{-(x_t,\lambda)}),
\end{equation*}
which completes the proof of
Theorem~\ref{thm:interv-multi}.\hfill$\Box$\bigskip

\paragraph{Proof of Lemma~\ref{lem:u-infty-diff}}
The differentiability of $u^\lambda_t$ with respect to $\lambda$ and
the ODE satisfied by its derivatives are classical results on the
regularity of the flow of ODEs (see e.g.\ Perko~\cite{Perko91}).

Moreover, since $\triangledown_\eta u^\lambda_t$ and $v^\lambda_t$ are
both solution of the ODE (\ref{eq:def-v}) (with different initial
conditions), it is trivial to transport the properties of
$v^\lambda_t$ proved in the proof of Theorem~\ref{thm:t-infty} to
$\triangledown_\eta u^\lambda_t$. In particular, $e^{-\mu
  t}\triangledown_\eta u^\lambda_t$ converges as $t\rightarrow+\infty$
to a non-zero vector $w^\lambda_\eta$ which is proportional to $\xi$.  We
only have to check that $\triangledown_\eta {\tilde
  u}^\lambda_\infty$ exists and that $w^\lambda_\eta=\triangledown_\eta {\tilde
  u}^\lambda_\infty$. Moreover, as for~(\ref{eq:ineq-v}),
\begin{equation*}
  \exp\Big(-\frac{c}{2}\int_0^t\sup_iu_{s,i}^\lambda\: ds\Big)e^{Dt}\eta \leq
  \triangledown_\eta
  u_t^\lambda\leq\exp\Big(-\frac{c}{2}\int_0^t\inf_iu_{s,i}^\lambda\:
  ds\Big)e^{Dt}\eta.  
\end{equation*}
Therefore, since $\exp(Dt)\geq 0$,
\begin{equation*}
  |\triangledown_\eta u^\lambda_t| \leq e^{Dt}|\eta|,
\end{equation*}
which implies that $e^{-\mu t}\triangledown_\eta u^\lambda_t$ is
uniformly bounded for $t\geq 0$, $\lambda\geq 0$ and $\eta$ in a
compact subset of $\mathbb{R}^k$.

Now, letting $t\rightarrow+\infty$ in~(\ref{eq:indep-CI-u-infty}) one
gets for any $h\geq 0$,
\begin{align*}
  |{\tilde u}_\infty^{\lambda+h\eta}-{\tilde u}_\infty^\lambda-P\int_0^h e^{-\mu
    t_0}\triangledown_\eta  u^{\lambda+r\eta}_{t_0} dr | 
  & =|{\tilde u}_\infty^{\lambda+h\eta}-{\tilde u}_\infty^\lambda-Pe^{-\mu
    t_0}(u_{t_0}^{\lambda+h\eta}-u_{t_0}^\lambda)| \\
  & \leq(1-e^{-c\varepsilon/2})P(e^{-\mu t_0}u^\lambda_{t_0}+e^{-\mu
    t_0}u^{\lambda+h\eta}_{t_0}) .
\end{align*}
Letting $\varepsilon\rightarrow 0$ (and thus $t_0\rightarrow+\infty$)
in the previous inequality, Lebesgue's convergence theorem yields
\begin{equation*}
  {\tilde u}_\infty^{\lambda+h\eta}-{\tilde u}_\infty^\lambda=\int_0^h
  w^{\lambda+r \eta}_\eta dr.
\end{equation*}
Therefore, $\tilde{u}^\lambda_\infty$ is differentiable with respect
to $\lambda$ and $\triangledown_\eta {\tilde
  u}^\lambda_\infty=w^\lambda_\eta$. The proof of
Lemma~\ref{lem:u-infty-diff} is completed.\hfill$\Box$\bigskip

\section{Some decomposable cases conditioned by different remote survivals}
\label{sec:particular-case}

In this section, we study some models for which the mutation matrix
$D$ is not irreducible: it is called `reducible' or `decomposable'.
In this case, the general theory developed above does not apply.  In
contrast with the irreducible case, the asymptotic behavior of the
MDW process and the MF diffusion depends on the type.\\
Decomposable critical multitype pure branching processes (without
motion and renormalization) were the subject of several works since
the seventies. See e.g.
\cite{Og75JMKU,FN76SSA,FN78/79ZWVG,Su81OJM,Zub,VaSa}.

\subsection{A first critical model}
\label{subsec:first model}

Our first example is a 2-types DW process with a reducible mutation
matrix of the form
\begin{equation}
  \label{eq:CP-D}
  D=\left(
    \begin{array}{rl}
      -\alpha & \alpha \\ 0 & 0
    \end{array}
  \right), \quad \alpha>0 .
\end{equation}
For this model type 1 (resp.\ type 2) is subcritical (resp.\
critical). Moreover mutations can occur from type 1 to type 2 but no
mutations from type 2 to type 1 are allowed.\\
In this section we analyze not only the law $\PP^*$ of MDW process
conditioned on the non-extinction of the whole population, but also
the MDW process conditioned on the survival of each type separately.
\begin{thm}
  \label{thm:cond-CP}
  Let $\PP$ be the distribution of the MDW process $X$ with mutation
  matrix~\textup{(\ref{eq:CP-D})} and non-zero initial condition $m$.
  Let us define $\PP^*, \hat{\PP}^*$ and $\check{\PP}^*$ for any $t>0$
  and $B\in{\cal F}_t$ by
  \begin{align*}
    \PP^*(B) & =\lim_{\theta\rightarrow\infty}\PP(B\mid \langle
    X_{t+\theta},{\bf 1}\rangle>0) \\
    \hat{\PP}^*(B) & =\lim_{\theta\rightarrow\infty}\PP(B\mid \langle
    X_{t+\theta,1},1 \rangle>0)
    \quad (\textrm{if }m_1 \not= 0)\\
    \check{\PP}^*(B) & =\lim_{\theta\rightarrow\infty}\PP(B\mid
    \langle X_{t+\theta,2}, 1 \rangle>0).
  \end{align*}
  Then, all these limits exist and, for any $t>0$,
  \begin{align}
    \check{\PP}^*\bigr|_{{\cal F}_t} = &  \PP^*\bigr|_{{\cal F}_t}
    = \frac{\langle X_t,\mathbf{1}\rangle}{\langle
      m,\mathbf{1}\rangle} \, \PP\bigr|_{{\cal F}_t} 
    \label{eq:cond-1} \\
    \textrm{ and }  \quad   \hat{\PP}^*\bigr|_{{\cal F}_t}
    = & \frac{\langle X_{t,1},1 \rangle}{\langle
      m_1,1 \rangle}\, e^{\alpha t} \, \PP\bigr|_{{\cal F}_t}
    . \label{eq:cond-2}
  \end{align}
\end{thm}

\paragraph{Proof}
Let us first prove~(\ref{eq:cond-2}). Using the method leading
to~(\ref{eq:bla}), we get that if $m_1\not=0$
$$
\hat{\PP}^*(B) =\lim_{\theta\rightarrow\infty}\EE\Big(\mathbb{1}_B
\frac{(x_t,\lim_{\lambda_1\rightarrow\infty}
  u_\theta^{(\lambda_1,0)})}
{(x,\lim_{\lambda_1\rightarrow\infty}
  u_{t+\theta}^{(\lambda_1,0)})}\Big).
$$
The cumulant $u_t^\lambda$ of the mass process satisfies for any
$\lambda = (\lambda_1,\lambda_2)$
\begin{equation*}
  \left\{
    \begin{array}{ll}
      \displaystyle{\frac{du_{t,1}^\lambda}{dt}=-\alpha
      u_{t,1}^\lambda+\alpha u_{t,2}^\lambda
        -\frac{c}{2}(u_{t,1}^\lambda)^2} & \quad u_{0,1}^\lambda=\lambda_1 \\
      \displaystyle{\frac{du_{t,2}^\lambda}{dt}=
        -\frac{c}{2}(u_{t,2}^\lambda)^2} & \quad u_{0,2}^\lambda=\lambda_2 .
    \end{array}
  \right.
\end{equation*}
The second equation admits as solution
\begin{equation}
  \label{eq:u_2}
  u_{t,2}^\lambda=\frac{\lambda_2}{1+\frac{c}{2}\lambda_2 t}
\end{equation}
and $u_{t,1}^\lambda$ can be computed explicitly if $\lambda_2=0$:
\begin{equation}
  \label{eq:u_1-lbd_2=0}
  u_{t,1}^{(\lambda_1,0)} =\frac{\lambda_1 e^{-\alpha
  t}}{1+\frac{c\lambda_1}{2\alpha}(1-e^{-\alpha t})} \quad \textrm{ and } \quad u_{t,2}^{(\lambda_1,0)}=0.
\end{equation}
We then get
$$
\hat{\PP}^*(B)= e^{\alpha t}\, \EE\big(\mathbb{1}_B \frac{x_{t,1}}{x_1} \big),
$$
which yields (\ref{eq:cond-2}).

Concerning $ \check{\PP}^*$, remark first that it is well defined even
if $m_2=0$ (but $m_1\not=0$) since particles of type 2 can be created
by particles of type 1.\\
We are going to prove~(\ref{eq:cond-1}) by a similar method as
Theorem~\ref{thm:cond}. Let us first compute $\xi$. The matrix $D$ has
two eigenvalues, $0$ and $-\alpha$, each of them with one-dimensional
eigenspace.  The (normalized) right and left eigenvectors of the
greatest eigenvalue $\mu=0$ are respectively
$\xi=(\frac{1}{2};\frac{1}{2})$ and $\eta=(0;2)$. Since $\xi>0$, the
proof of Lemma~\ref{lem:pte-u} (and therefore Lemma~\ref{lem:pte-u}
itself) is still valid for this specific matrix, except for the
assertion (i), which has to be reduced to the following: if
$\lambda_2>0$, then, for any $t>0$, $u_t^\lambda>0$ (if $\lambda_2=0$,
then $u^\lambda_{t,2}\equiv 0$).

Therefore, as in the proof of Theorem~\ref{thm:cond}, we can prove that
\begin{equation*}
  \forall \lambda\geq 0,\ \forall t\geq 0,\quad
  |u_t^\lambda-e^{Dt}\lambda|\leq K\|\lambda\|te^{Dt}\lambda 
\end{equation*}
and thus that, for $\theta$ and $t_0$ such that
$\|u_{t_0}^\lambda\|\leq 1/K(\theta+t)$,
\begin{multline}
  \left|\frac{( a,u_{t_0+\theta}^\lambda)}
    {( b,u_{t_0+\theta+t}^\lambda)}
    -\frac{( a,e^{D\theta}u_{t_0}^\lambda)}
    {( b,e^{D(\theta+t)}u_{t_0}^\lambda)}\right| \\
  \leq
  \frac{2K\|a\|\|u_{t_0}^\lambda\|\theta\|e^{D\theta}u_{t_0}^\lambda\|}
  {( b,e^{D(\theta+t)}u_{t_0}^\lambda)}
  +\frac{2K\|a\|\|b\|\|e^{D\theta}u_{t_0}^\lambda\|
    \|u_{t_0}^\lambda\|(t+\theta)\|e^{D(\theta+t)}u_{t_0}^\lambda\|}
  {( b,e^{D(\theta+t)}u_{t_0}^\lambda)^2}.
\label{eq:et-une-autre-eq}
\end{multline}
Let us compute explicitly the exponential of the matrix $Dt$.
Since $D^n=(-\alpha)^n N$ where $ N=\left(
  \begin{array}{lr} 1 & -1 \\ 0 & 0 \end{array}\right)$, then
\begin{equation*}
  e^{Dt}=P+e^{-\alpha t}N,
  \mbox{\ with\ }
  P=(\xi_i\eta_j)_{1\leq
    i,j\leq 2}=\left(
    \begin{array}{ll}
      0 & 1 \\ 0 & 1
    \end{array}\right).
\end{equation*}
This has the same form as~(\ref{equivalentexpDt}) in the irreducible
case, except that $P\not > 0$. Because of this, we cannot obtain a
bound for~(\ref{eq:et-une-autre-eq}) uniform in $\lambda\in\RR^2$ as
in the proof of Theorem~\ref{thm:cond}. However, we can restrict to a
subset of $\RR^2$ for which the convergence is uniform and which
covers the two limits involved in the computation of $\check{\PP}^*$
($\lambda_1=0$ and $\lambda_2\rightarrow+\infty$) and $\PP^*$
($\lambda_1\rightarrow+\infty$ and $\lambda_2\rightarrow+\infty$).

This can be done as follows: if $\lambda_1=\lambda_2$, then
$u_{t,1}^\lambda = u_{t,2}^\lambda$ for any $t\geq 0$. Therefore, for
any $\lambda\not=\{0\}$ such that $\lambda_2\geq\lambda_1$,
$u_{t,2}^\lambda\geq u_{t,1}^\lambda>0$ for any $t>0$ (if at some time
these quantities are equal, they remain equal for any larger time).
Then, since
$$
e^{D\theta}u_{t_0}^\lambda =
(u_{t_0,2}^\lambda + e^{-\alpha\theta} (u_{t_0,1}^\lambda -
u_{t_0,2}^\lambda), u_{t_0,2}^\lambda),
$$
this quantity converges when $\theta\rightarrow \infty$ to
$(u_{t_0,2}^\lambda,u_{t_0,2}^\lambda)$, \emph{uniformly} in $\lambda$
such that $\lambda_2>\lambda_1\geq 0$.

From this follows as in the proof of Theorem~\ref{thm:cond} that
\begin{equation*}
  \lim_{\theta\rightarrow\infty}
  \frac{u_\theta^\lambda}{(x,u_{t+\theta}^\lambda)}
  =\frac{\xi}{(x,\xi)}
\end{equation*}
uniformly for $\lambda$ in the set of $(\lambda_1 ; \lambda_2)\not=0$
such that $\lambda_2\geq\lambda_1\geq 0$. This ends the proof
of~(\ref{eq:cond-1}).\hfill$\Box$\bigskip

From the local density of $\PP^*$ (resp.  $\hat{\PP}^*$) with respect
to $\PP$, we easily obtain as in Section~\ref{sec:immigration} the
following expressions for the Laplace functionals of the different
conditioned processes.

\begin{thm}
  \label{thm:V-CP}
  The probability measure $\PP^*(=\check{\PP}^*)$ is characterized by
  \begin{equation*}
    \forall f \in  C_b(\RR,\RR^2)_+  ,\quad \EE^*(\exp-\langle
    X_t,f\rangle\mid X_0=m)
    =\frac{\langle m,V_tf\rangle}{\langle m,\xi\rangle
    }e^{-\langle m,U_tf\rangle}
  \end{equation*}
  where $V_tf$ is the unique semigroup solution of the PDE
  \textup{(\ref{eq:def-V})}.
  The probability measure $\hat{\PP}^*$ is characterized by
  \begin{equation*}
    \forall f \in  C_b(\RR,\RR^2)_+  ,\quad \hat{\EE}^*(\exp-\langle
    X_t,f\rangle\mid X_0=m)
    =\frac{\langle m,\hat{V}_tf\rangle}{\langle m_1,1\rangle}
    e^{\alpha t}e^{-\langle m,U_tf\rangle}
  \end{equation*}
  where $\hat{V}_tf$ satisfies the same PDE as $V_tf$ except for the
  initial condition $ \hat{V}_0f=(1 ; 0)$.
\end{thm}

\subsection{Long time behaviors of the Feller diffusions}
\label{sec:cumul-CP}

Let us now analyze the long time behavior of the various conditioned
MF diffusions.
\begin{prop}
  \label{prop:CP}
  \begin{description}
  \item[\textmd{(a)}] The first type vanishes in $\PP^*$-probability \
    when $t\rightarrow\infty$, that is
    \begin{equation*}
      \forall \varepsilon>0, \quad \lim_{t\rightarrow\infty} \PP^*(
      x_{t,1} > \varepsilon) =0  .
    \end{equation*}
  \item[\textmd{(b)}] Under $\hat{\PP}^*$, the first type
    converges in distribution to the probability measure $\Gamma(2,
    2\alpha/c)$
    \begin{equation*}
      \lim_{t\rightarrow\infty}\hat{\PP}^*(x_{t,1}\in\cdot)
      \overset{(d)}{=}\Gamma(2,2\alpha/c).
    \end{equation*}
  \item[\textmd{(c)}] The second type $x_{t,2}$ explodes
    in $\PP^*$- and in $\hat{\PP}^*$-probability when
    $t\rightarrow\infty$.
  \end{description}
\end{prop}

\paragraph{Proof}
We first compute the vector $v_t^\lambda:=V_t\lambda,
\lambda\in\RR^2$.
\begin{equation*}
  \label{eq:TL-v-CP}
  \EE^*(\exp-( x_t,\lambda)\mid x_0=x)
  =\frac{( x,v_t^\lambda)}{(
    x,\xi)}e^{-( x,u_t^\lambda)}
\end{equation*}
with
\begin{equation*}
  \left\{
    \begin{array}{ll}
      \displaystyle{\frac{dv_{t,1}^\lambda}{dt}
        =-\alpha v_{t,1}^\lambda+\alpha v_{t,2}^\lambda
        -c u_{t,1}^\lambda v_{t,1}^\lambda,} & \quad v_{0,1}^\lambda=1/2, \\
      \displaystyle{\frac{dv_{t,2}^\lambda}{dt}=-c u_{t,2}^\lambda
        v_{t,2}^\lambda,} & \quad v_{0,2}^\lambda=1/2.
    \end{array}\right.
\end{equation*}
Therefore, replacing $u_{t,2}^\lambda$ by its value obtained
in~(\ref{eq:u_2}),
\begin{equation*}
  v_{t,2}^\lambda=\frac{1}{2(1+\frac{c}{2}\lambda_2 t)^2}
\end{equation*}
and
\begin{equation}
  \label{eq:v_1-CP}
  v_{t,1}^\lambda=\frac{1}{2}e^{-\alpha t-c\int_0^tu_{s,1}^\lambda ds}
  \Big(1+\alpha\int_0^t 
  \frac{e^{\alpha s+c\int_0^s u_{\tau,1}^\lambda
      d\tau}}{(1+\frac{c}{2}\lambda_2 s)^2} \,  ds \Big).
\end{equation}
In particular, if $\lambda_2=0$, $v_{t,2}^{(\lambda_1; 0)}=1/2$ and
one gets from the explicit expression~(\ref{eq:u_1-lbd_2=0}) of
$u_{t,1}^{(\lambda_1;0)}$
\begin{equation*}
  v_{t,1}^{(\lambda_1,0)}=\frac{1}{2}-\frac{c\lambda_1
    (1+\frac{c}{2\alpha}\lambda_1)te^{-\alpha t}
    +\frac{c^2}{\alpha^2}\lambda_1^2 e^{-2\alpha t}}
  {2\left(1+\frac{c}{2\alpha}\lambda_1 (1-e^{-\alpha t})\right)^2}.
\end{equation*}
Similarly, for $\hat{v}_t^\lambda:=\hat{V}_t\lambda, \lambda\in\RR^2$,
we get
\begin{equation*}
  \label{eq:TL-hat-v}
  \hat{\EE}^*(\exp-( x_t,\lambda)\mid x_0=x)
  =\frac{( x,\hat{v}_t^\lambda)}{( x,\xi)}e^{\alpha
    t}e^{-( x,u_t^\lambda)}
\end{equation*}
with
\begin{equation*}
  \label{eq:def-hat-v}
  \left\{
    \begin{array}{ll}
      \displaystyle{\frac{d\hat{v}_{t,1}^\lambda}{dt}
        =-\alpha\hat{v}_{t,1}^\lambda+\alpha\hat{v}_{t,2}^\lambda
        -c u_{t,1}^\lambda \hat{v}_{t,1}^\lambda,}
      & \quad \hat{v}_{0,1}^\lambda=1, \\
      \displaystyle{\frac{d\hat{v}_{t,2}^\lambda}{dt}=-c u_{t,2}^\lambda
        \hat{v}_{t,2}^\lambda,} & \quad \hat{v}_{0,2}^\lambda=0.
    \end{array}\right.
\end{equation*}
Therefore, $\hat{v}_{t,2}^\lambda=0$ and
\begin{equation} \label{eq:hat-v-1-general} \hat{v}_{t,1}^\lambda=\exp
  \Big(-\alpha t-c\int_0^t u_{s,1}^\lambda \, ds \Big) .
\end{equation}
In particular, if $\lambda_2=0$, using~(\ref{eq:u_1-lbd_2=0}) again,
\begin{equation}
  \label{eq:hat-v-1}
  \hat{v}_{t,1}^{(\lambda_1;0)} =\frac{e^{-\alpha t}}
  {\left(1+\frac{c}{2\alpha}\lambda_1 (1-e^{-\alpha t})\right)^2}.
\end{equation}

Now (a) and (b) can be deduced from the facts that $\lim_{t\rightarrow
  \infty}v_t^{(\lambda_1;0)}=\xi$ and $\lim_{t\rightarrow
  \infty}e^{\alpha t}\hat{v}_t^{(\lambda_1;0)}
=(1/(1+c\lambda_1/2\alpha)^2;0)$.\\
The explosion of $x_{t,2}$ in $\PP^*$-probability in~(c) is a
consequence of the fact that $\lim_{t\rightarrow \infty}e^{\alpha
  t}\hat{v}_t^{(\lambda_1;\lambda_2)}=
\big(\exp( -c\int_0^\infty u_{s,1}^\lambda \, ds );0 \big)=(0;0)$, by Lemma~\ref{lem:pte-u} (iii). \\
Finally, it follows from~(\ref{eq:v_1-CP}) that
\begin{align*}
  v_{t,1}^\lambda & \leq\frac{e^{-\alpha t}}{2}+\frac{\alpha}{2}\int_0^t
  \frac{e^{-\alpha(t-s)-c\int_s^tu_{\tau,1}^\lambda d\tau}}
  {(1+\frac{c}{2}\lambda_2 s)^2}ds \\
  & \leq\frac{e^{-\alpha t}}{2}+\frac{e^{-\alpha t/2}}{2}
  +\frac{\alpha}{c\lambda_2(1+\frac{c}{4}\lambda_2 t)}
\end{align*}
where the last inequality is obtained by splitting the integral over
the time interval $[0,t]$ into the sum of the integrals over
$[0,\frac{t}{2}]$ and $[\frac{t}{2},t]$.  This implies the first part
of~(c).\hfill$\Box$\bigskip

We interpret this proposition as follows. Conditionally on the
survival of the whole population, the weakest type gets extinct and
the strongest type has the same behavior as in the critical monotype
case. Conversely, conditionally on the long time survival of the
weakest type, the weakest type behaves at large time as in the
monotype subcritical case and the strongest type explodes.

\subsection{A  more general subcritical decomposable model }
\label{sec:CP-more-gal}

We consider a generalization of the previous model. The mutation matrix is now given by
\begin{equation}  \label{eq:CP-D+}
D=\left(
    \begin{array}{rr}
      -\alpha & \alpha \\ 0 & -\beta
    \end{array}
  \right)
\end{equation}
where $\alpha>0$ (as before) and $\beta>0$ with $\beta \not= \alpha $.\\

In this case, the whole population is subcritical. Here again,
mutations are only possible from type 1 to type 2. If $\beta<\alpha$,
type 2 is ``less subcritical'' than type 1 (as in the
previous case) but if $\alpha<\beta$, type 1 is ``less subcritical'' than
type 2. We will see below that the behavior of the various conditioned
processes is strongly related to the so-called dominating type, which
is the first one if $\alpha<\beta$ and the second type if
$\beta<\alpha$.\\
Before treating separately both cases with different techniques, we
define the common ingredients we need.

We can easily compute the normalized right eigenvector $\xi$ for the
greatest eigenvalue $\mu$. If $\beta<\alpha$, $\mu=-\beta$ and
$\xi=\frac{1}{2\alpha-\beta}(\alpha;\alpha-\beta)$ and if
$\alpha<\beta$, $\mu=-\alpha$ and $\xi=(1;0)$. We can also explicitly
compute the exponential of the mutation matrix:
\begin{equation*}
  e^{Dt}=e^{-\beta t}\left(
    \begin{array}{lc}
      0 & \frac{\alpha}{\alpha-\beta} \\ 0 & 1
    \end{array}\right)+e^{-\alpha t}\left(
    \begin{array}{lc}
      1 & -\frac{\alpha}{\alpha-\beta} \\ 0 & 0
    \end{array}\right).
\end{equation*}
The cumulant $u_t^\lambda$ of the mass process satisfies
\begin{equation} \label{eq:u-t-1}
  \left\{
    \begin{array}{ll}
      \displaystyle{\frac{du_{t,1}^\lambda}{dt}=
        -\alpha u_{t,1}^\lambda + \alpha u_{t,2}^\lambda
        -\frac{c}{2}(u_{t,1}^\lambda)^2} & \quad u_{0,1}^\lambda=\lambda_1 \\
      \displaystyle{\frac{du_{t,2}^\lambda}{dt}=
        -\beta u_{t,2}^\lambda-\frac{c}{2}(u_{t,2}^\lambda)^2} & \quad
      u_{0,2}^\lambda=\lambda_2 
    \end{array}
  \right.
\end{equation}
Thus $u_{t,2}^\lambda$ is given by
\begin{equation}
  \label{eq:expl-u_2}
  u_{t,2}^\lambda=\frac{\lambda_2e^{-\beta
      t}}{1+\frac{c}{2\beta}\lambda_2 (1-e^{-\beta t})} .
\end{equation}
One can compute $u_{t,1}^\lambda$ explicitly only when $\lambda_2=0$,
and in this case, as in ~(\ref{eq:u_1-lbd_2=0}),
 
\begin{equation*}
  u_{t,1}^{(\lambda_1;0)}=\frac{\lambda_1 e^{-\alpha
      t}}{1+\frac{c}{2\alpha}\lambda_1 (1-e^{-\alpha t})} \quad
  \textrm{ and } \quad u_{t,2}^{(\lambda_1;0)}=0.
\end{equation*}

We now consider the system of equations 
\begin{equation}
  \label{eq:y-v}
  \left\{
    \begin{array}{l}
      \displaystyle{\frac{dh_{t,1}}{dt}
        =-\alpha h_{t,1}+\alpha h_{t,2}-cu_{t,1}^\lambda h_{t,1}} \\
      \displaystyle{\frac{dh_{t,2}}{dt}=-\beta h_{t,2}-c u_{t,2}^\lambda
        h_{t,2}}
    \end{array}\right.
\end{equation}
which solutions are given by
\begin{equation}
  \label{eq:y_2}
  h_{t,2}=\frac{h_{0,2}\, e^{-\beta t}}{\big(
    1+\frac{c\lambda_2}{2\beta}(1-e^{-\beta t}) \big)^2}.
\end{equation}
and
\begin{equation}
  \label{eq:y_1}
  h_{t,1} = e^{-\alpha t-c\int_0^tu_{s,1}^\lambda ds}
  \Big(h_{0,1}+\alpha\int_0^te^{\alpha s+c\int_0^s u_{\tau,1}^\lambda d\tau}
  h_{s,2}ds \Big).
\end{equation}

We denote as before by $v^\lambda_t$, $\hat{v}^\lambda_t$ or
$\check{v}^\lambda_t$ the respective solutions of (\ref{eq:y-v}) with
initial conditions $v^\lambda_0=\xi$, $\hat{v}^\lambda_0=(1;0)$ and
$\check{v}^\lambda_0=(0;1)$.

\subsubsection{Case $\beta<\alpha$}
\label{sec:beta<alpha}

We now identify the laws obtained by conditioning with respect to the
various remote survivals.
 
\begin{thm}
  \label{cond-CP-2}
  Let $\PP^*$ (resp. $ \hat{\PP}^*, \check{\PP}^*$) be the conditioned
  laws defined in Theorem~\textup{\ref{thm:cond-CP}}
  where $\PP$ is the law of the MDW process with mutation matrix
  given by~\textup{(\ref{eq:CP-D+})} with $\beta<\alpha$ and non-zero
  initial condition $m$. It holds
  \begin{align*}
    \check{\PP}^*\bigr|_{{\cal F}_t}=\PP^*\bigr|_{{\cal F}_t}
    = &\frac{\langle X_t,\xi\rangle}{\langle
      m,\xi\rangle}e^{\beta t} \, \PP\bigr|_{{\cal F}_t} \\
    \textrm{and } \quad \hat{\PP}^*\bigr|_{{\cal F}_t}
    =& \frac{\langle X_{t,1},1 \rangle}{\langle
      m_1,1 \rangle}e^{\alpha t} \, \PP\bigr|_{{\cal F}_t} \quad
    (\textrm{if }m_1 \not= 0).
  \end{align*}
\end{thm}

\paragraph{Sketch of the proof}
The greatest eigenvalue of $D$ is $\mu=-\beta$, the normalized right
eigenvector for $\mu $ is
$\xi=\frac{1}{2\alpha-\beta}(\alpha;\alpha-\beta)$ and the normalized
left eigenvector is
$\eta=(0;\frac{2\alpha-\beta}{\alpha-\beta})$.\\
As in the proof of Theorem~\ref{thm:cond-CP}, $\xi>0$, so that
Lemma~\ref{lem:pte-u} holds (except assertion~(i)\:) and we can use a
similar method. The only difficulty is to find a domain $E\subset
\RR^2_+$ such that, for each initial condition $\lambda \in E$, the
cumulant semigroup $u_t^\lambda$ takes its values in $E$ and
$\{\lambda_1/\lambda_2, \lambda \in E\}$ is bounded.  To this aim, one
can check that, if $0\leq
u_{t,1}^\lambda=\frac{\alpha}{\alpha-\beta}u_{t,2}^\lambda$ at some
time $t\geq 0$, then $\frac{du^\lambda_{t,1}}{dt}
\leq\frac{\alpha}{\alpha-\beta} \frac{du^\lambda_{t,2}}{dt}$.
Therefore, if $0\leq\lambda_1\leq\frac{\alpha}{\alpha-\beta}\lambda_2$
with $\lambda_2>0$, one has $0 < u^{\lambda}_{t,1} \leq
\frac{\alpha}{\alpha-\beta} u^\lambda_{t,2}$ for any positive
$t$.\hfill$\Box$\bigskip

Let us now analyze the behavior for large $t$ of the mass process
under the three measures $\PP^*$, $ \hat{\PP}^*$ and $ \check{\PP}^*$.
Since $\hat{v}_{0,2}=0$, $\hat{v}_{t,2}\equiv 0$ as in
Section~\ref{sec:cumul-CP} and~(\ref{eq:hat-v-1-general}) holds.
Therefore the behavior of $x_t$ under $\hat{\PP}^*$ is exactly the
same as for $\beta=0$, treated in Proposition~\ref{prop:CP}~(b) and~(c).\\
The long time behavior of $x_t$ under $\PP^*$ is different from
Section~\ref{sec:cumul-CP} and is given in the following proposition.
\begin{prop}
  \label{prop:CP-bis}
  \begin{description}
  \item[\textmd{(a)}] The first type vanishes in $\PP^*$-probability
    when $t\rightarrow\infty$.
  \item[\textmd{(b)}] Under $\PP^*$,  the second type
    converges in distribution to the probability measure $\Gamma(2,
    2\beta/c)$
    \begin{equation*}
      \lim_{t\rightarrow\infty}\PP^*(x_{t,2}\in\cdot)
      \overset{(d)}{=}\Gamma(2,2\beta/c).
    \end{equation*}
  \end{description}
\end{prop}

\paragraph{Proof}
Since
\begin{equation*}
  \EE^*(e^{-( x_t,\lambda)} \mid x_0=x)
  =\frac{( x,v_t^\lambda)}{( x,\xi)}
e^{\beta t}e^{-( x,u_t^\lambda)}
\end{equation*}
we have to compute $\lim_{t\rightarrow\infty} v_t^{\lambda}e^{\beta t}$.\\
For the proof of (a) we remark that, from (\ref{eq:y_2}) and
(\ref{eq:y_1}), $ v_{t,2}^{(\lambda_1;0)}=\xi_2\, e^{-\beta t}$ and
\begin{equation*}
  v_{t,1}^{(\lambda_1;0)} = e^{-\alpha t}\exp(-c\int_0^t u_{s,1}^{(\lambda_1;0)} ds)
  \Big(\xi_1 + \alpha \, \xi_2 \int_0^t e^{(\alpha-\beta) s}
  \exp(c\int_0^s u_{\tau,1}^{(\lambda_1;0)} d\tau) ds \Big).
\end{equation*}
Since 
$$
\exp(-c \int_0^t u_{s,1}^{(\lambda_1;0)} ds)
=\exp\Big(-c \int_0^t \frac{\lambda_1 e^{-\alpha
    s}}{1+\frac{c}{2\alpha}\lambda_1 (1-e^{-\alpha s})} ds\Big)
= \frac{1}{\big(1+\frac{c}{2\alpha}\lambda_1 (1-e^{-\alpha t})\big)^2}
$$
one obtains 
\begin{align*} 
  v_{t,1}^{(\lambda_1;0)} e^{\beta t}
  & = \exp\Big(-\frac{(\alpha - \beta) t}
  {(1+\frac{c}{2\alpha}\lambda_1 (1-e^{-\alpha t}))^2}\Big) \xi_1 \\
  &  + \frac{\alpha  }{\big(1+\frac{c}{2\alpha}\lambda_1 (1-e^{-\alpha t})\big)^2}
  \int_0^t e^{-(\alpha - \beta)(t-s)}
  \big(1+\frac{c}{2\alpha}\lambda_1 (1-e^{-\alpha s})\big)^2
  ds \, \xi_2 .
\end{align*} 
The integral can be computed explicitly and is equal, for $t$ large,
to
$$
\frac{1}{\alpha - \beta} (1+\frac{c}{2\alpha}\lambda_1 )^2 +
O(e^{-(\alpha - \beta)t}) \, .
$$   
Thus, 
$$
\lim_{t\rightarrow\infty} v_{t,1}^{(\lambda_1;0)} e^{\beta t} =
\frac{\alpha}{\alpha - \beta} \, \xi_2 = \xi_1 .
$$
For the proof of (b), it suffices to show that  
$$
\lim_{t\rightarrow \infty} v_t^{(0;\lambda_2)}e^{\beta t}
=\frac{1}{(1+\frac{c}{2\beta}\lambda_2)^2} \,  \xi .
$$
From (\ref{eq:y_2}), it is clear that $\lim_{t\rightarrow \infty}
v_{t,2}^{(0;\lambda_2)}\, e^{\beta
  t}=\frac{1}{(1+\frac{c}{2\beta}\lambda_2)^2} \xi_2$. It then
remains to compute the limit, for $\lambda_2 >0$, of
$v_{t,1}^{(0;\lambda_2)}$ as $t\rightarrow\infty$.
Using~(\ref{eq:y_1}), we get
\begin{equation*}
  v^{(\lambda_1;\lambda_2)}_{t,1} e^{\beta t} =e^{-(\alpha-\beta)t-c\int_0^tu^\lambda_{s,1}ds} \, \xi_1
  + \alpha\xi_2\int_0^t\frac{e^{-(\alpha-\beta)(t-s)-c{\int_s^tu^\lambda_{\tau,1}d\tau}}}
  {\big(1+\frac{c\lambda_2}{2\beta}(1-e^{-\beta s})\big)^2} ds.
\end{equation*}
The first term is $O(e^{-(\alpha - \beta)t})$ and goes to 0 as
$t\rightarrow\infty$. The limit of the integral can be computed as
follows:
\begin{multline*}
  \Big|\int_0^t\frac{e^{-(\alpha-\beta)(t-s)-c{\int_s^tu^\lambda_{\tau,1}d\tau}}}
  {\big( 1+\frac{c}{2\beta}\lambda_2 (1-e^{-\beta s})\big)^2}ds
  -\frac{1}{(1+\frac{c}{2\beta}\lambda_2)^2}
  \int_0^te^{-(\alpha-\beta)(t-s)}ds\Big| \\
  \begin{aligned}
    & \leq \bar K\int_0^te^{-(\alpha-\beta)(t-s)}
    \Big|1-e^{-c\int_s^tu^\lambda_{\tau,1}d\tau}
    \Big(1-\frac{\frac{c\lambda_2}{2\beta}e^{-\beta s}}
    {1+\frac{c\lambda_2}{2\beta}}\Big)^{-2}\Big|ds \\
    & \leq \bar K \bigg(e^{-(\alpha-\beta)t/2}
    +\frac{t}{2}\Big(1-e^{-c\int_{t/2}^{+\infty}u^\lambda_{s,1}ds}\Big)
    \vee\Big(\Big(1-\frac{\frac{c\lambda_2}{2\beta}e^{-\beta
        t/2}}{1+\frac{c}{2\beta}\lambda_2}\Big)^{-2}-1\Big) \bigg)
  \end{aligned}
\end{multline*}
where the positive constant $\bar K$ may vary from line to line and
where the last inequality is obtained by splitting the integration
over the time intervals $[0,\frac{t}{2}]$ and $[\frac{t}{2},t]$.\\
Now, by Lemma~\ref{lem:pte-u} (ii), $\lim_{t\rightarrow \infty}
\int_{t/2}^{+\infty} u^\lambda_{s,1}ds =0$.  Therefore,
$$
\lim_{t\rightarrow\infty} v^{\lambda}_{t,1} e^{\beta t} 
= \frac{\alpha\xi_2}{\big(1+\frac{c}{2\beta}\lambda_2 \big)^2}
\int_0^\infty e^{-(\alpha-\beta)s}ds 
=\frac{\xi_1}{\big(1+\frac{c }{2\beta}\lambda_2 \big)^2}
$$
as required.\hfill$\Box$\bigskip

Here again, this result can be interpreted as follows: for $i=1,2$,
conditionally on the survival of the type $i$, this type $i$ behaves
as if it was alone, and the other type $j$ explodes or goes extinct
according to whether it is stronger or weaker.

\subsubsection{Case $\alpha<\beta$}
\label{sec:alpha<beta}

When $\alpha < \beta$, the greatest eigenvalue of $D$ is $\mu=-\alpha$
and the normalized right eigenvector to $\mu $ is $\xi=(1;0)$. In
particular, $\xi\not >0$, so that Lemma~\ref{lem:pte-u} does not hold
and we cannot use the previous method anymore. However, in our
specific example, $u^\lambda_{t,2}$ can be explicitly computed, and,
by~(\ref{eq:u-t-1}), $u^\lambda_{t,1}$ is solution of the
one-dimensional differential equation
\begin{equation}
  \label{eq:u-gal}
  \frac{dy_t}{dt}=-\alpha y_t-\frac{c}{2}y_t^2
  +\frac{\alpha \lambda_2 e^{-\beta t}}{1+\frac{c}{2\beta} \lambda_2
    (1-e^{-\beta t})} .
\end{equation}
This equation can be (formally) extended to the case
$\lambda_2=\infty$ as
\begin{equation}
  \label{eq:u-gal-lbd-infty}
  \frac{dy_t}{dt}=-\alpha y_t-\frac{c}{2}y_t^2
  +\frac{2\alpha\beta e^{-\beta t}}{c (1-e^{-\beta t})} .
\end{equation}
The following technical lemma gives (non-explicit) long-time estimates
of the solutions of~(\ref{eq:u-gal}) that are sufficient to compute
the various conditioned laws of the MDW process. We postpone its proof
at the end of the subsection.

\begin{lem}
  \label{lem:last-pty-u}
  For any $\lambda_2\in[0,\infty]$, let ${\cal Y}(\lambda_2)$ denote
  the set of solutions $y_t$ of~\textup{(\ref{eq:u-gal})} (or
  of~\textup{(\ref{eq:u-gal-lbd-infty})} if $\lambda_2=\infty$)
  defined (at least) on $(0,\infty)$. For any $y\in{\cal
    Y}(\lambda_2)$, the limit $C(\lambda_2,y):=
  \lim_{t\rightarrow\infty} y_t e^{\alpha t} $ exists and satisfies
  $$
  0< \inf_{\lambda_2 \ge 1,\: y\in{\cal Y}(\lambda_2)} C(\lambda_2,y) \le 
  \sup_{\lambda_2 \ge 1,\: y\in{\cal Y}(\lambda_2)} C(\lambda_2,y) < + \infty
  $$
\end{lem}

We now identify the laws obtained by conditioning $\PP$ with respect
to the various remote survivals.
 
\begin{thm}
  \label{thm:cond-last-CP}
  Let $\PP^*$ (resp. $ \hat{\PP}^*, \check{\PP}^*$) be the conditioned
  laws defined in Theorem~\textup{\ref{thm:cond-CP}}, where $\PP$ is the
  MDW process with mutation matrix given by \textup{(\ref{eq:CP-D+})}
  with $\alpha<\beta$ and initial condition $m=(m_1;m_2)$ with
  $m_1\not=0$. It holds
  \begin{equation*}
    \hat{\PP}^*\bigr|_{{\cal F}_t}=\check{\PP}^*\bigr|_{{\cal F}_t}
    =\PP^*\bigr|_{{\cal F}_t}
    =\frac{\langle X_{t},\xi \rangle}{\langle
      m,\xi \rangle} \, e^{\alpha t} \, \PP\bigr|_{{\cal F}_t} .
  \end{equation*}
  When $m_1=0$ and $m_2\not=0$, 
  \begin{equation*}
    \check{\PP}^*\bigr|_{{\cal F}_t}
    =\PP^*\bigr|_{{\cal F}_t}
    =\frac{\langle X_{t,2},1 \rangle}{\langle
      m_2, 1\rangle} \, e^{\beta t} \, \PP\bigr|_{{\cal F}_t} .
  \end{equation*}
  and $\hat{\PP}^*$ is not defined.
\end{thm}

\paragraph{Proof}
Our usual method consists in computing the following limits as
$\theta$ goes to $+\infty$ :
\begin{equation}
  \label{eq:3-limits}
  \frac{\lim_{\lambda_1\rightarrow\infty}u_\theta^{(\lambda_1;0)}}
  {(x,\lim_{\lambda_1\rightarrow\infty}u_{t+\theta}^{(\lambda_1;0)})},\quad
  \frac{\lim_{\lambda_2\rightarrow\infty}u_\theta^{(0;\lambda_2)}}
  {(x,\lim_{\lambda_2\rightarrow\infty}u_{t+\theta}^{(0;\lambda_2)})},\quad
  \frac{\lim_{\lambda_1,\lambda_2\rightarrow\infty}u_\theta^{(\lambda_1;\lambda_2)}}
  {(x,\lim_{\lambda_1,\lambda_2\rightarrow\infty}u_{t+\theta}^{(\lambda_1;\lambda_2)})}.
\end{equation}
It is elementary to prove that, as monotone limits of solutions
of~(\ref{eq:u-gal}), the function $ \theta \mapsto
\lim_{\lambda_1\rightarrow\infty}u_{\theta,1}^{(\lambda_1;0)}$ is
still solution of~(\ref{eq:u-gal}) with $\lambda_2=0$, and the
functions $ \theta \mapsto \lim_{\lambda_2\rightarrow\infty}
u_{\theta,1}^{(0;\lambda_2)} $ and $ \theta \mapsto
\lim_{\lambda_1,\lambda_2\rightarrow\infty}u_{\theta,1}^{(\lambda_1;\lambda_2)}$
are solutions of~(\ref{eq:u-gal-lbd-infty}) (a priori only defined for
$t>0$).\\
Therefore, we can use Lemma~\ref{lem:last-pty-u} and the explicit
formula~(\ref{eq:expl-u_2}) for $u^\lambda_{t,2}$ to compute the three
limits of~(\ref{eq:3-limits}). In each case, the dominant term is the
one including $u^\lambda_{t,1}$, except when $m_1=0$, where the only
remaining term is the one including $u^\lambda_{t,2}$.\hfill$\Box$\bigskip

Finally, we give the long time behavior of the mass process under
$\PP^*$ (which is equal to $\check{\PP}^*$ and $\hat{\PP}^*$ when
this last measure exists).

\begin{prop}
  \label{prop:last-long-time} 
  \begin{description}
  \item[\textmd{(a)}] If $m_1\not=0$, the laws under $\PP^*$ of the
    mass of both types $x_{t,1}$ and $x_{t,2}$ converge when
    $t\rightarrow \infty$. More precisely
    \begin{equation*}
      \lim_{t\rightarrow\infty}\PP^*(x_{t,1}\in\cdot)
      \overset{(d)}{=}\Gamma(2,2\alpha/c).
    \end{equation*}
    $x_{t,2}$ converges in $\PP^*$-distribution to a non-trivial (and
    non-explicit) distribution on $\RR_+$.
  \item[\textmd{(b)}] If $m_1=0$ ($m_2\not=0$), $x_{t,1}\equiv 0\ \
    \PP^*-a.s.$ and
    \begin{equation*}
      \lim_{t\rightarrow\infty}\PP^*(x_{t,2}\in\cdot)
      \overset{(d)}{=}\Gamma(2,2\beta/c).
    \end{equation*}
  \end{description}
\end{prop}

\paragraph{Proof}
With the previous notation, when $m_1\not=0$,
\begin{equation*}
  \EE^*(\exp-( x_t,\lambda))
  =\frac{( x,\hat{v}_t^\lambda)}{( x,\hat{v}^\lambda_0)}e^{\alpha
    t}e^{-( x,u_t^\lambda)}.
\end{equation*}
Since $\hat{v}_{t,2}=0$ and $\hat{v}_{t,1}=\exp(-\alpha
t-c\int_0^tu^\lambda_{s,1}ds)$,
\begin{equation*}
  \lim_{t\rightarrow\infty}e^{\alpha t}\hat{v}_t=\Big(\exp-\int_0^\infty
  u^\lambda_{s,1}ds;0\Big)
\end{equation*}
where $\exp-\int_0^\infty u^\lambda_{s,1}ds \in (0,1)$ by
Lemma~\ref{lem:last-pty-u}. In order to prove the convergence in
distribution of $x_{t,2}$, it remains to prove that
$$\lim_{\lambda_2\rightarrow 0}\int_0^\infty
u^{(0,\lambda_2)}_{s,1}ds=0.$$
Because of~(\ref{eq:u-gal}), $\dot{u}^\lambda_{t,1}\leq -\alpha
u^\lambda_{t,1}+\alpha\lambda_2e^{-\beta t}$. Therefore, it is easy to
check that $u^{\lambda}_{t,1}\leq (\beta+2)\lambda_2e^{-\beta t}$ for
all $t\geq 0$ if $\lambda_1\leq (\beta+2)\lambda_2$ (simply
differentiate the difference). This implies the required result.

When $\lambda_2=0$, the computations can be made explicitly as in the
proof of Proposition~\ref{prop:CP} (b) and give the usual Gamma limit
distribution for $x_{t,1}$ under $\PP^*$ when $t\rightarrow+\infty$.

If $m_1=0$ (i.e.\ $x_1=\langle m_1,1\rangle =0$),
\begin{equation*}
  \EE^*(\exp-( x_t,\lambda))
  =\frac{( x,\check{v}_t^\lambda)}{( x,\check{v}^\lambda_0)}e^{\beta
    t}e^{-( x,u_t^\lambda)}=
  \check{v}_{t,2}^\lambda e^{\beta
    t}e^{-( x,u_t^\lambda)}
\end{equation*}
and the computations are the same as in the monotype case.\hfill$\Box$\bigskip

We interpret this last result as follows: when the first type is
present, it dominates the asymptotic behavior of both types, since its
subcriticality is weaker than the one of the second type, although
mutations from type 2 to type 1 do not occur.

\begin{rem}
  \label{rem:case-alpha=beta}
  \textup{Theorem~\ref{thm:cond-last-CP} and
    Proposition~\ref{prop:last-long-time} are still valid in the remaining case
    $\alpha=\beta$. In this case, it is actually possible to be more precise than in
    Lemma~\ref{lem:last-pty-u} by proving, using a similar method,
    that a solution $y_t$ to~(\ref{eq:u-gal})
    or~(\ref{eq:u-gal-lbd-infty}) (with $\alpha=\beta$) satisfies
    $y_t\sim C(\lambda_2)te^{-\alpha t}$ where
    $C(\lambda_2)=\alpha\lambda_2/(1+c\lambda_2/2\alpha)$ if
    $\lambda_2<\infty$ and $C(\infty)=2\alpha^2/c$.}
\end{rem}

\paragraph{Proof of lemma~\ref{lem:last-pty-u}}
Let $z_t:=e^{\alpha t}y_t$. It solves the equation
\begin{equation}
  \label{eq:z}
  \frac{dz_t}{dt}=-\frac{c}{2}e^{-\alpha t}z_t^2
  +\frac{\alpha\lambda_2e^{-(\beta-\alpha)t}}
  {1 + \frac{c}{2\beta} \lambda_2 (1-e^{-\beta t})} .
\end{equation}

Let us first prove that $z_t$ is bounded for $t\in[1,\infty)$,
uniformly in $\lambda_2$ and independently of the choice of the
solution $y_t$ of~(\ref{eq:u-gal}) or~(\ref{eq:u-gal-lbd-infty})
defined on $(0,\infty)$. For any $t\geq 1$ and
$\lambda_2\in[0,\infty]$,
\begin{equation*}
  \frac{dz_t}{dt}\leq\frac{2\alpha\beta}{c(1-e^{-\beta})}e^{-(\beta-\alpha)t}.
\end{equation*}
Since the integral of the above r.h.s.\ over $[1,\infty)$ is finite,
we only have to prove that $z_1$ is bounded uniformly in
$\lambda_2\geq 0$ and independently of the choice of $y_t$.  Now, for
any $t\in[\frac{1}{2},1]$ and $\lambda_2\in[0,\infty]$,
\begin{equation*}
  \frac{dz_t}{dt}\leq\frac{2\alpha\beta}{c}
  \frac{e^{-(\beta-\alpha) t}}{1-e^{-\beta /2}}-\frac{ce^{-\alpha}}{2}z_t^2.
\end{equation*}
In particular, for any $t\in[\frac{1}{2},1]$, the first term in the
r.h.s.\ above is bounded and bounded away from 0. Thus, there exists a
constant $K$ such that, if $z_t\geq K$ and $\frac{1}{2} \leq t\leq 1$,
$\frac{dz_t}{dt} \leq -\frac{ce^{-\alpha}}{4}z_t^2$. Therefore,
distinguishing between $z_{1/2}\leq K$ and $z_{1/2}>K$, one obtains
\begin{equation*}
  z_1\leq\frac{z_{1/2}}{1+\frac{c}{4}e^{-\alpha}z_{1/2}}\vee
  K\leq\frac{4e^\alpha}{c}\vee K=:K'<\infty.
\end{equation*}

Second, it follows from~(\ref{eq:z}) and from the boundedness of $z_t$
that $|dz_t/dt|\leq K''(e^{-(\beta-\alpha)t}+e^{-\alpha t})$ for some
constant $K''$ for any $t\geq 1$. Therefore, $z_t$ converges as
$t\rightarrow\infty$. Moreover, this function is uniformly bounded
from above for $t\geq 1$ by some constant $K'''$ independent of the
particular function $z_t$ considered. Therefore, since $dz_t/dt\geq
-cK'''e^{-\alpha t}z_t/2$ for $t\geq 1$, the limit of $z_t$ when
$t\rightarrow +\infty$ is also greater than
$z_1\exp(-cK'''e^{-\alpha}/2)$. 

Then, it only remains to prove that $z_1$ is bounded away from 0,
uniformly in $\lambda_2\in[1,\infty]$.  This follows from the fact
that, for any $\lambda_2\geq 1$, there exists a constant $M>0$
independent of $\lambda_2$ and $t$ such that, for
$t\in[\frac{1}{2},1]$,
\begin{equation*}
  \frac{dz_t}{dt}\geq M-\frac{c}{2}z_t^2.
\end{equation*}
This implies that there exists $M'$ such that, if $z_t\leq M'$ for 
$t\in[\frac{1}{2},1]$, $dz_t/dt\geq M/2$, and thus
\begin{equation*}
  z_1\geq \big( z_{1/2}+\frac{M}{4} \big)\wedge M'\geq
  \frac{M}{4}\wedge M'>0,
\end{equation*}
which completes the proof of Lemma~\ref{lem:last-pty-u}.\hfill$\Box$\bigskip

\subsection{Exchange of long time limits}
\label{sec:interv-CP}

As in the irreducible case, one can study the interchangeability of
the long time limit ($t\rightarrow+\infty$) of the conditioned Feller diffusion 
and the limit of long time survival
($\theta\rightarrow+\infty$).  The same method as in
Proposition~\ref{thm:interv-multi} yields, for $i,j\in\{1,2\}$ and
$\lambda \in \RR_+$,
\begin{equation*}
  \lim_{t\rightarrow\infty}\EE(e^{-\lambda x_{t,i}}\mid
  x_{t+\theta,j}>0)=\lim_{t\rightarrow\infty}\frac{(
    x,u_t^{\lambda^i + \lim_{\bar\lambda\rightarrow\infty}u_{\theta}^{\bar\lambda^j }}
    -u_t^{\lambda^i})}{(
    x,\lim_{\bar\lambda\rightarrow\infty}u_{t+\theta}^{\bar\lambda^j})}
\end{equation*}
where $\lambda^1=(\lambda;0)$ and $\lambda^2=(0;\lambda)$.

However, the computation of these quantities requires precise
information about the behavior of $u_t^\lambda$ as a function of its
initial condition $\lambda$ for $t$ large. The cases we could handle
are the one with non degenerate limits.  In the model introduced in
section \ref{subsec:first model}, it corresponds to $i=j=1$ and the
computation reduces to the monotype case studied in
Proposition~\ref{prop:interv-CP-monotype}. In the model introduced in
section \ref{sec:CP-more-gal}, with $\beta < \alpha$, it corresponds
to $i=j=1$ and $i=j=2$, and with $\alpha < \beta$, to all cases. For
$i=j=1$ the proof is based on explicit expressions like in the
monotype case, and for the other cases, the arguments are similar to those
used in the proof of Theorem~\ref{thm:interv-multi} (except in the
case $\alpha < \beta$ and $m_1=0$, where the computation can
also be done explicitly). To summarize:

\begin{prop}
  \label{prop:interv-CP}
  In the cases described above, one can interchange both limits in time:
  \begin{align*}
    \lim_{\theta\rightarrow\infty}\lim_{t\rightarrow\infty}
    \PP(x_{t,1}\in\cdot \mid x_{t+\theta,1}>0)
    & \overset{(d)}{=} \lim_{t\rightarrow\infty}\lim_{\theta\rightarrow\infty}
    \PP(x_{t,1}\in\cdot \mid x_{t+\theta,1}>0) \\
    & \overset{(d)}{=}\Gamma(2,2\alpha/c) \quad \textrm{ if  } \, 0 \leq \beta<
    \alpha,\\
    \lim_{\theta\rightarrow\infty}\lim_{t\rightarrow\infty}
    \PP(x_{t,2}\in\cdot \mid x_{t+\theta,2}>0) & \overset{(d)}{=}
    \lim_{t\rightarrow\infty}\lim_{\theta\rightarrow\infty}
    \PP(x_{t,2}\in\cdot \mid x_{t+\theta,2}>0) \\
    & \overset{(d)}{=}\Gamma(2,2\beta/c) \quad \textrm{ if  } \, 0 < \beta< \alpha\\
    & \phantom{\overset{(d)}{=}\Gamma(2,2\beta/c) \quad} \textrm{ or }\,
    0<\alpha < \beta \textrm{ and }\, m_1=0,
  \end{align*}
  \begin{align*}
    \lim_{\theta\rightarrow\infty}\lim_{t\rightarrow\infty}
    \PP(x_{t,1}\in\cdot \mid x_{t+\theta,1}>0) & \overset{(d)}{=}
    \lim_{\theta\rightarrow\infty}\lim_{t\rightarrow\infty}
    \PP(x_{t,1}\in\cdot \mid x_{t+\theta,2}>0)\\ \overset{(d)}{=}
    \lim_{t\rightarrow\infty}\lim_{\theta\rightarrow\infty}
    \PP(x_{t,1}\in\cdot \mid x_{t+\theta,1}>0) & \overset{(d)}{=}
    \lim_{t\rightarrow\infty}\lim_{\theta\rightarrow\infty}
    \PP(x_{t,1}\in\cdot \mid x_{t+\theta,2}>0) \\
    & \overset{(d)}{=}\Gamma(2,2\alpha/c) \quad \textrm{ if  } \, 0 < \alpha < \beta
  \end{align*}
  and
  \begin{align*}
    \lim_{\theta\rightarrow\infty}\lim_{t\rightarrow\infty}
    \PP(x_{t,2}\in\cdot \mid x_{t+\theta,1}>0) & \overset{(d)}{=}
    \lim_{\theta\rightarrow\infty}\lim_{t\rightarrow\infty}
    \PP(x_{t,2}\in\cdot \mid x_{t+\theta,2}>0)\\ \overset{(d)}{=}
    \lim_{t\rightarrow\infty}\lim_{\theta\rightarrow\infty}
    \PP(x_{t,2}\in\cdot \mid x_{t+\theta,1}>0) & \overset{(d)}{=}
    \lim_{t\rightarrow\infty}\lim_{\theta\rightarrow\infty}
    \PP(x_{t,2}\in\cdot \mid x_{t+\theta,2}>0)
  \end{align*}
  if $0 < \alpha < \beta$ and $m_1\not =0$ (in this last case the
  limit is not known explicitly).
\end{prop}

\paragraph{Acknowledgments}
The first author is grateful to the DFG, which supported his Post-Doc
in the Dutch-German Bilateral Research Group ``Mathematics of Random
Spatial Models from Physics and Biology'', at the Weierstrass
Institute for Applied Analysis and Stochastics in Berlin, where part
of this research was made.

\end{document}